\documentclass[a4paper, 10pt]{amsart}

\usepackage{amsmath, latexsym,  amssymb, amscd}
\usepackage[all]{xy}
\newcommand{\ecuno}{\epsilon_1(C,\eta)}
\newcommand{\ecdue}{\epsilon_2(C,\eta)}

\newcommand{\centro}{{\rm{A}}}
\newcommand{\finito}{{\rm{B}}}

\newcommand{\expaltc}{{1+\frac{1}{2(2g-3)}}}

\newcommand{\expm}{{4n}}
\newcommand{\expalt}{{1+\frac{1}{2n}}}
\newcommand{\expaltu}{{1+\frac{1}{2n}}}
\newcommand{\gain}{{\frac{1}{2n}}}
\newcommand{\valeta}{{2^4(g+s)(g-1)^2}}
\newcommand{\coeta}{{8(g+s)(g-1)}}

\newcommand{\cqp}{{(C(\overline{\qe})\times \gamma)}}
\newcommand{\cp}{{C\times \gamma}}

\newcommand{\pippo}{\psi}
\newcommand{\fiffo}{\phi}

\newcommand{\enne}{{b}}
\newcommand{\ti}{{f}}

\newcommand{\cunop}{{c_p}}
\newcommand{\ezerop}{{\varepsilon_p}}
\newcommand{\czero}{{c_0}}

\newcommand{\cuno}{c_1}

\newcommand{\cdue}{c_2}

\newcommand{\zono}{{\rm {End}}(E) }
\newcommand{\emor}{{\rm {End}} }
\newcommand{\ring}{R}

\newcommand{\cod}{{\rm {cod}} }
\newcommand{\g}{\theta}

\newcommand{\qe}{\mathbb{Q}}
\newcommand{\co}{\mathbb{C}}

\newcommand{\ze}{\mathbb{Z}}

\newcommand{\upxi}{\zeta}

\newcommand{\punto}{y}

\newcommand{\boundaC}{\min\left(1, \frac{\kzero}{g},\left(\frac{\epsilon(C)}{g\kzero}\right)^{8g^2}\right)}
\newcommand{\boundb}{\min\left(\euno,\frac{K_1}{g},\left(\frac{\epsilon(C)}{g\kuno}\right)^{8(g+s)g}\right)}

\newcommand{\boundctra}{\frac{1}{g^{4g}}\min(1,K_0^{-1})^{4g}\min\left(1,\kzero, \left(\frac{\epsilon(C)}{g\kzero}\right)^{8g^2}\right)^{4g}}

\newcommand{\ezero}{\varepsilon_0(p,\tau)}
\newcommand{\euno}{\varepsilon_1}

\newcommand{\edue}{{\varepsilon_2}}
\newcommand{\eduep}{{\varepsilon_3}}

\newcommand{\etre}{\varepsilon'_0(p,\tau)}

\newcommand{\ecinque}{\varepsilon_4}

\newcommand{\kzero}{K_0}

\newcommand{\kuno}{K_1}

\newcommand{\kdue}{K_2}
\newcommand{\kduep}{K_3}

\newcommand{\ktre}{K_4}

\newtheorem{thm}{Theorem}[section]
\newtheorem{con}[thm]{Conjecture}
\newtheorem{propo}[thm]{Proposition}
\newtheorem{lem}[thm]{Lemma}

\newtheorem{D}[thm]{Definition}
\newtheorem{thm1}[thm]{Theorem}

\newtheorem*{propa*}{Proposition B}
\newtheorem*{propb*}{Proposition A}

\newtheorem{remark}[thm]{Remark}

\title[ {The intersection of a curve with a union of translated 
 subgroups }]
{ The intersection of a curve with a union of translated codimension $2$
subgroups in a power of an elliptic curve}      

\author[{ Evelina Viada}]{  
 }

\begin{document}

\maketitle
\centerline{Evelina Viada\footnote{Evelina Viada,
     Universit\'e de Friboug Suisse, P\'erolles, D\'epartement de Math\'ematiques, 23, Chemin du Mus\'ee, CH-1700 Fribourg,
Suisse, viada@math.ethz.ch, 
    evelina.viada@unifr.ch.}
    \footnote{Supported by the SNF (Swiss National Science Foundation).}
}
\begin{abstract}
Let $E$ be an elliptic curve. Consider
an irreducible algebraic curve $C$ embedded in $E^g$.  The curve  is transverse if it is not contained in any
translate of a proper algebraic subgroup of $E^g$. Furthermore $C$ is
weak-transverse if it is not contained in  any proper algebraic
subgroup. Suppose that both $E$ and $C$  are  defined over
the algebraic numbers. 

We prove that the algebraic points of a transverse curve $C$ which are close to the union of all algebraic subgroups of $E^g$ of codimension $2$ translated by points in a subgroup $\Gamma$ of $E^g$  of finite rank are a set of bounded height. The notion of close is defined using a height function. 
If $\Gamma$ is trivial, it is sufficient to suppose that $C$ is weak-transverse. 

Then, we introduce a method to determine the finiteness of these sets.
From a conjectural lower bound for the normalised height of a transverse curve $C$, we    deduce that the  above sets  are finite. At present, such a lower bound exists for $g\le3$.

Our  results are  optimal, for what concerns the codimension of the algebraic subgroups. 
\end{abstract}

\section{introduction}

We present the problems in a general context.

Denote by  $A$  a semi-abelian variety   over $\overline{\qe}$  of
dimension $g$. Consider an irreducible algebraic  subvariety  $V$
of  $A$, defined over $\overline{\qe}$. 
We say that 
\begin{itemize}
 \item $V$ is transverse if $V$ is not contained in any
translate of a proper algebraic subgroup of $A$.

\item  $V$ is
weak-transverse if $V$ is not contained in  any proper algebraic
subgroup  of $A$.

\end{itemize}

Given an integer $ r$ with $1 \le r \le g$ and a subset $F$ of $A(\overline{\qe})$, we define the set
$$S_r(V,F)=V(\overline{\qe}) \cap \left(\bigcup_{\mathrm{cod}B \ge r} B+F\right),$$
where $B$ varies over all semi-abelian subvarieties of $A$ of
codimension at least $r$ and $$B+F=\{b+f \,\,\,: \,\,\,b\in B, \,\,\,f\in F\}.$$ 
For $r>g$, we define $S_r(V,F)$ to be the empty set. We  denote the set $S_r(V, A_{\rm {Tor}})$
simply  by $S_r(V)$.
Note that $$S_{r+1}(V,F)\subset S_{r}(V,F).$$

A natural question to ask would be;  for which sets $F$ and integers
$r$, the set 
 $S_r(V,F)$ is non Zariski-dense in $V$.
 
Sets of this kind, for $r=g$, appear in the literature in the
context of the Mordell-Lang, of the Manin-Mumford and of the Bogomolov Conjectures.
More recently  Bombieri, Masser and Zannier \cite{BMZ}  have proven that the set $S_2(C)$ is finite,
for a transverse
 curve $C$ in a torus. They
investigate,  for the first time,  intersections with the union of all
algebraic subgroups of a given codimension. This opens a vast
number of conjectures for subvarieties of semi-abelian varieties.

In this article we consider the elliptic case for curves.  
Let $E$ be an elliptic curve  and $C$ an irreducible algebraic curve
 in $E^g$, both defined over $\overline{\qe}$.  
 Let $||\cdot||$ be a semi-norm on $E^g(\overline{\qe})$ induced by a height function.
For $\varepsilon \ge 0$, we denote $$\mathcal{O}_{\varepsilon}=\{ \xi \in E^g(\overline{\qe}) : ||\xi|| \le \varepsilon\}.$$ 
Let $\Gamma\subseteq E^g(\overline{\qe})$ be a subgroup of finite rank.  Define
$\Gamma_{\varepsilon}=\Gamma + \mathcal{O}_\varepsilon.$

\begin{con}
\label{geni}
Let $C\subset E^g$. Then,

\begin{enumerate}

\item If $C$ is weak-transverse,  $S_2(C)$  is finite.

\item If $C$ is transverse,  $S_2(C,\Gamma)$ is finite.

\item  If $C$ is weak-transverse,  there exists  $\varepsilon>0$ such that $S_2(C,\mathcal{O}_\varepsilon)$ is finite.

\item If $C$ is transverse,  there exists  $\varepsilon>0$ such that $S_2(C,\Gamma_\varepsilon)$ is finite.

\end{enumerate}
\end{con}

The strong
hypotheses of $C$ transverse or the weak hypotheses of $C$ weak-transverse is a crucial difference in the setting. Please take note as to which hypotheses is assumed in the different statements.

Clearly  iv. implies  ii. by setting $\varepsilon=0$, and similarly  iii. implies  i.

 The
union of all algebraic subgroups of codimension $g$ is exactly
 the torsion of $E^g$. Then, $C\cap \Gamma_\varepsilon\subset S_g(C, \Gamma_\varepsilon)\subset S_2(C,\Gamma_\varepsilon)$.  So, Conjecture \ref{geni} 
iii. and  iv. imply  the Bogomolov  (\cite{ulm}, \cite{zang})
and  the Mordell-Lang plus Bogomolov (\cite{poonen}) Theorems  respectively.

Partial results related to  i. and  ii. have been proven.
In \cite{V} we solve a weak form of  i. There we assume the stronger hypothesis that  $C$ is  transverse. If $E$
has C.M. then $S_2(C)$ is finite. If $E$ has no  C.M. then
$S_{\frac{g}{2}+2}(C)$ is finite. In \cite{RV}  R\'emond and the
author present a weak version of   ii. Again if $E$ has C.M. the
result is optimal. If $E$ has no C.M. the codimension of the algebraic
subgroups  depends on $\Gamma$.  In   addition,  we show that  i. and  ii. are equivalent.
 Note that  there are no trivial implications between  iii. and  iv.,  because of the different hypotheses on $C$.

These  known  proofs rely on Northcott's theorem: a set is finite if and only if  it has bounded height and degree. To prove that the degree is bounded one  uses  Siegel's Lemma and    an essentially optimal generalized Lehmer's Conjecture.  Up to a logarithm factor, the generalized Lehmer conjecture is  presently known for a point in  a torus \cite{fra1} and in a C.M. abelian
 variety \cite{sinh}.  This method has  some disadvantages: first it is only known  to work  for transverse curves and for $\varepsilon=0$, secondly 
  a quasi optimal  generalized Lehmer's Conjecture   is not likely to be proven in a near future for a general abelian variety. 

In this article we introduce a different method. First, we 
bound the height also for weak-transverse curves. 
\begin{thm}
\label{altzero}
There exists  $\varepsilon>0$
  such that:

\begin{enumerate}
\item  If $C$ is weak-transverse,  $S_2(C,\mathcal{O}_\varepsilon)$ has bounded height.

\item If $C$ is transverse, $S_2(C,\Gamma_\varepsilon)$ has bounded height.

\end{enumerate}
\end{thm}
 
 The proof of both statements  uses a Vojta inequality as stated in
\cite{RV}  Proposition 2.1.
 The second assertion is proven in \cite{RV} Theorem 1.5. To prove
 the first assertion (see  section \ref{bound}), we embed $S_2(C,\mathcal{O}_\varepsilon)$ into two sets
associated to a transverse curve. We then manage to  apply a Vojta inequality on each of these two sets. 

As a second result, we prove:

\begin{thm}
\label{equivalenza}
For $r\ge 2$, the following statements are equivalent:
\begin{enumerate}

\item  If $C$ is weak-transverse, then there exists  $\varepsilon>0$ such that $S_r(C,\mathcal{O}_\varepsilon)$ is finite.

\item If $C$ is transverse, then there exists  $\varepsilon>0$ such that $S_r(C,\Gamma_\varepsilon)$ is finite.

\end{enumerate}
\end{thm}
This theorem is not as easy as  the equivalence of  Conjecture \ref{geni} part i. and ii. That i. implies ii. is quite elementary. 
The other implication is
delicate. In particular we make use of Theorem \ref{altzero} (see section \ref{bound}).

In the third instance, we
show how to avoid the use of the Siegel Lemma and of the generalized Lehmer Conjecture.
Instead, we use Dirichlet's Theorem and a conjectural effective
version of the Bogomolov Theorem.
Bogomolov's Theorem  states that  the set of points of small height on a curve of genus at least 2 is finite.
We define  $\mu(C)$ as the supremum
of the   reals $\epsilon(C)$ such that $S_g(C,\mathcal{O}_{\epsilon(C)})=C\cap
\mathcal{O}_{\epsilon(C)}$ is finite. The essential minimum of $C$ is  $\mu(C)^2$ (note that in the literature, often, the notation $\mathcal{O}_\varepsilon$ corresponds to the set, we denote in this work, $\mathcal{O}_{\varepsilon^2}$; thus in the references given below the bounds are given for the essential minimum and not for its square root $\mu(C)$ as we use here).

 Non-optimal effective lower bounds
for $\mu(C)$  are given by S. David and P. Philippon  \cite{sipa}  Theorem 1.4 and \cite{spec} Theorem 1.6.
The lower bound we need is  the elliptic analogue  to a theorem of Amoroso and David.  In \cite{fra} Theorem 1.4, they prove an essentially optimal lower bound for a variety in a torus.  The following conjecture is a weak form of \cite{spec} Conjecture 1.5 part ii. where the line bundle is fixed.

\begin{con}
\label{sin2}
Let $A=E_1\times \cdots \times E_g$ be a product of elliptic curves defined over a number field $k$. Let $L$ be the tensor product of the pull-backs of symmetric line bundles on $E_i$ via the natural projections.
Let $C\subset A$ be an irreducible transverse curve defined over $\overline{\qe}$.  
Let $\eta$ be any  positive real.
Then, there exists a constant $ c(g,A,\eta)=c(g,\deg_{{{L}}} A, h_L(A), [k:\mathbb{Q}],\eta) $   such that, 
for  $$\epsilon(C,\eta) = {c(g,A,\eta)} {(
\deg_{{{L}}} C)^{-\frac{1}{2(g-1)}-\eta}},$$
the set 
$$C(\overline{\mathbb{Q}}) \cap \mathcal{O}_{\epsilon(C,\eta)}$$ is finite.

\end{con}

In section \ref{cg}, we prove:

\begin{thm}
\label{tutto}
Conjecture \ref{sin2}  implies Conjecture \ref{geni}.
\end{thm}

Conjecture \ref{sin2} can be stated for subvarieties of $A$. Galateau \cite{aur} proves that such a cunjecture holds for varieties of codimension $1$ or $2$ in a product of elliptic curves. Then, for $g\le 3$, Conjecture \ref{geni} holds unconditionally. 

Theorems \ref{altzero} and \ref{tutto} are  optimal with respect to the codimension of the algebraic subgroups - see remark \ref{optimal}.

We have already pointed out that Conjecture \ref{geni}  implies  the Bogomolov Conjecture and the Mordell-Lang plus Bogomolov Theorem. Let us emphasise that our Theorem \ref{tutto} does not give a new proof of the Bogomolov Conjecture, as we assume such an effective result. On the other hand,  it gives a new proof of the Mordell-Lang plus Bogomolov Theorem, under the assumption of Conjecture \ref{sin2}.

The  proof  of Theorem \ref{tutto} is based on two steps.
A union of  infinitely many sets is finite if and only if 
\begin{itemize}
 \item[(1)] the  union can be taken over    finitely
many sets,
\item[(2)] all sets in the union are
finite.
\end{itemize}

(1)  is a 
typical problem of Diophantine approximation. The proof relies on  Dirichlet's Theorem on the rational approximation of reals. 
The fact that we consider small neighbourhoods enables us to move  the algebraic subgroups `a bit'. So we can   consider only subgroups of bounded degree, which are finitely many (see Proposition $\centro$, \S \ref{subunion}).

The second step (2) places itself in the context of  the  height theory. Its proof relies on Conjecture \ref{sin2}. 
The bound  $\epsilon(C,\eta)$ depends on the invariants of the ambient variety and  on the degree of $C$.
 A weaker dependence  on the degree of $C$ would not be enough for our application. Also the non-dependence of the bound on the  field of definition of  $C$ proves useful. Playing on Conjecture \ref{sin2}, we produce a sharp lower bound for the essential minimum of the image of a curve under certain morphisms (see Proposition $\finito$, \S \ref{due}).

The effectiveness aspect of our method is  noteworthy; the  use of  a Vojta inequality makes
Theorem \ref{altzero}, and consequently  Theorem \ref{tutto},
ineffective. 
Though, the rest of the method   is effective.  Indeed, in section \ref{effettivo},
we  prove a weaker, but effective analogue of
Theorem \ref{tutto}.

\begin{thm}
\label{Ctrao}
Assume Conjecture \ref{sin2}. Let $C$ be transverse. Then, there exists an effective $\varepsilon>0$ such that  the set  $S_2(C,\mathcal{O}_\varepsilon)$ is finite.

\end{thm}

 A bound   for the number of
 points of small height on the curve would then imply a bound for the
 cardinality of   $S_2(C,\mathcal{O}_\varepsilon)$
 for $C$ transverse  and $\varepsilon$ small (see Theorem \ref{ord1}).

The toric version of Theorem \ref{Ctrao}  is independently  studied  by P. Habegger in his Ph.D. thesis  \cite{PHIL}.
He follows the idea of using a Bogomolov type bound, proven in the toric case in \cite{fra} Theorem 1.4.   He proves the finiteness of
$S_2(C,\mathcal{O}_{\epsilon})$, for $\epsilon>0$ and $C$ a transverse curve in a
torus.


 \newpage

{\it Acknowledgments} 

My  deepest thanks go to the Referee for his precious comments, which encouraged me  to drastically improve a first manuscript. 

I kindly thank Francesco Amoroso and \"Ozlem Imamoglu for  inspiring  me to pursue research in
Mathematics.

I thank the Centro de Giorgi and in particular Umberto Zannier for
inviting me  to the Diophantine Geometry Semester in 2005. I gave
there a talk where the core material of this paper was presented.

I thank  the SNF for the  financial support which  allows me to work serenely near my family.

It is to my husband Robert that I am most thankful, without him I would have given up long ago and I would not know how nice it is to be happy.

 \newpage

\section{preliminaries}

\label{nota}

\subsection{Morphisms and their height}
Let $(\ring, |\cdot|)$ be a hermitian ring, that means $R$ is a domain and $|\cdot|$ an absolute value on $R$.

We denote by $M_{r,g}(\ring)$ the  module of $r\times g$ matrices with
entries in  $\ring$.

For $F=(f_{ij})\in M_{r,g}(\ring)$,  we define the height of $F$ as the maximum of the absolute value of its entries
$$H(F)=\max_{ij}|f_{ij}|.$$

Let $E$ be an elliptic curve defined over a number field.
The ring of endomorphism $\emor(E)$ is isomorphic  either to $\ze$ (if $E$ does not have C.M.) or to an order in an imaginary quadratic field  (if $E$ has C.M.). 
 We consider on $\emor(E)$ the  standard absolute value of $\mathbb{C}$. Note that this  absolute value does not depend on the   embedding of $\emor(E)$ in $\mathbb{C}$.
 
We identify a morphism $\phi:E^g \to E^r$ with a matrix  in $M_{r,g}({\rm{End}}(E))$. Note that the set of morphism of height bounded by a constant is finite.

An  intrinsic definition of  absolute value on $\emor(E)$ can be given using  the Rosati-involution. 

In the following, we aim to be as transparent as possible, polishing statements from technicality. Therefore,
 we principally present proofs for $E$ without C.M. Then $\emor(E)$ is identified with $\ze$ and a morphism $\phi$ with an integral matrix. 
 In the final section, we  explain  how to deal with the technical
 complication of a ring of endomorphisms of rank $2$ and with a
 product of elliptic curves instead of a power. 

\subsection{Small points}
On $E$, we fix a symmetric very ample line bundle $\mathcal{L}$. On
$E^g$, we consider the bundle $L$ which is  the tensor product of the
pull-backs of $\mathcal{L}$ via the natural projections on the
factors. Degrees are computed with respect to the polarization $L$. 

Usually $E^g(\overline{\qe})$ is endowed with the $L$-canonical N\'eron-Tate height
$h'$. Though, to simplify constants, we prefer to 
define  on $E^g$ the height of the maximum 
\begin{equation*}
h(x_1,\dots ,x_g)=\max_i(h(x_i)).
\end{equation*}
where $h(\cdot)$  on $E(\overline{\qe})$  is  the  $\mathcal{L}$-canonical N\'eron-Tate height.
The height $h$ is the square of a norm  $||\cdot||$ on $E^g(\overline\qe)\otimes \mathbb{R}$.  For a point $x \in E^g(\overline\qe)$, we write $||x||$ for $||x\otimes1||$. 

Note that  $h(x)\le h'(x)\le gh(x)$. Hence, the two norms induced by $h$ and $h'$ are equivalent.

For $a\in \emor(E)$, we denote by $[a]$ the multiplication by $a$.
For  $y\in E^g(\overline{\qe})$ it holds  $$\Big|\Big|[a]y\Big|\Big|=|a|\cdot||y||.$$

The height  of a non-empty set $S\subset E^g(\overline\qe)$ is the supremum of the heights  of
its elements. The norm of $S$ is  the non-negative square root
of its height.

For $\varepsilon \ge 0$, we denote 
$$\mathcal{O}_{\varepsilon}=\mathcal{O}_{\varepsilon,E^g}=\{ \xi \in E^g(\overline{\qe}) : ||\xi|| \le \varepsilon\}.$$

\subsection{Subgroups}
Let $M$ be a  $\ring$-module. 
The $\ring$-rank of $M$ is the supremum  of the cardinality  of a set of $\ring$-linearly independent elements of $M$. 
If $M$ has finite rank $s$, a maximal free set of $M$ is a set of $s$ 
linearly independent  elements of $M $. 
If $M$  is a free  $\ring$-module of  rank  $s$,  
we call  a set of $s$ generators of $M$, integral generators of $M$.

Please note that a free $\ze$-module of finite rank is a lattice; in the literature, what we call integral generators can be called basis, and what we define as maximal free set is a basis of the vector space given by tensor product with the quotient field of $\ring$.

We say that $(M,||\cdot||)$ is a hermitian $R$-module if $M$ is an $R$-module and $||\cdot||$ is a norm on the tensor product of $M$ with the quotient field of $\ring$. For an element $p \in M$ we write $||p||$ for $||p\otimes 1||$.

Let $E$ be an elliptic curve. In the following, we will simply say module for an ${\rm{End}}(E)$-module.

Note that any subgroup of $E^g(\overline{\qe})$ of finite rank  is contained in a sub-module of  finite rank. Conversely,  a sub-module  of $E^g$ of finite rank is a subgroup of finite rank. 

Let $\Gamma$ be a subgroup of finite rank of $E^g(\overline{\qe})$.
We define
$$\Gamma_{\varepsilon}=\Gamma + \mathcal{O}_\varepsilon.$$
The  saturated module  ${\Gamma}_0$ of the coordinates group of $\Gamma$ (in short of $\Gamma$) is a sub-module of $E(\overline{\qe})$ defined as 
\begin{equation}
\label{gammazero}
{\Gamma}_0=\{\phi(y)  \in E {\rm {\,\,\,for\,\,\,}} \phi: E^g\to E {\rm{ \,\,\,\and \,\,\,}}Ny\in \Gamma{\rm \,\,\,with\,\,\,}  N\in \ze^*\}.
\end{equation}
Note that $\Gamma^{g}_0=\Gamma_0\times \dots \times \Gamma_0$ is a sub-module of $E^g$ invariant via the image
or preimage of isogenies. Furthermore, it contains $\Gamma$ and it is a module of finite rank. This shows that  to prove finiteness statements for $\Gamma$ it is enough to prove them for $\Gamma_0^{g}$. 

We denote by $s$ the rank of $\Gamma_0$. Let $\gamma_{1}, \dots ,\gamma_{s}$ be a maximal free set of $\Gamma_0$.
We denote the associated point of $E^s$ by
\begin{equation*}
\gamma=(\gamma_1,\dots, \gamma_{s}).
\end{equation*}

\section{ Some Geometry of Numbers}
We present a property of geometry of numbers, we then extend it to points of  $E^g(\overline{\qe})$. The  idea is that,  if  in $\mathbb{R}^n$  we
consider $n$ linearly independent vectors, and we
move them within a `small' angle, they will still be linearly
independent. Furthermore, the norm of a linear combination of such vectors 
depends on the norm of these vectors, on their angles, and on the
norm of the coefficients of the combination.

Such  estimates are frequent in the geometry of numbers.


The following lemma is a reformulation of  \cite{sc} Theorem 1.1 of Schlickewei or  \cite{V} Lemma 3.

\begin{lem}
\label{bg} 
Every  hermitian free $\mathbb{Z}$-module of rank $n$ admits   integral generators  $\rho_1, \dots , \rho_n$ such that 
for all integers $\alpha_i$ 
 
 \begin{equation*}
 {{\czero}}(n)\sum_i|\alpha_i|^2||\rho_i||^2\le \big|\big|\sum_i\alpha_i\rho_i\big|\big|^2
 \end{equation*}
with ${{\czero}}(n)$ a constant depending only on $n$.

\end{lem}
\begin{proof}
A hermitian free $\mathbb{Z}$-module $(\Gamma, ||\cdot||)$ of rank $n$ is a lattice in the metric space $\Gamma_{\mathbb{R}}$ given by tensor product  with $\mathbb{R}$. The proof is now equal to the proof of  \cite{V} Lemma 3 page 57, from line 19 onwards, where one shall read $n$ for $r$ and $\rho_i$ for $g_i$.
\end{proof}

This lemma allows us to explicit the comparison constant for two norms on a finite dimensional vector space over the quotient field of $R$.

\begin{propo}
\label{angoli}
Let $(M,||\cdot||)$ be a hermitian $R$-module, with $R$ a  finitely genereted free $\ze$-module.
Let $p_1, \dots ,p_s$ be  $R$-linearly independent elements of $M$. Then
there exists an effective positive constant $\cuno(p,\tau)$    such that, 
for all  $b_1,\dots,b_s \in R$, it holds
$$\cuno(p,\tau)\sum_i|b_i|_R^2||p_i||^2 \le \big|\big|\sum_ib_ip_i\big|\big|^2,$$
where $p=(p_1, \dots ,p_s)$ and $\tau=(1,\tau_2, \dots, \tau_t)$  integral generators of $R$.

\end{propo}
\begin{proof}

The sub-module of $M$ defined by $\Gamma_{\mathbb{Z}}=\langle p_1, \dots, p_s, \dots, \tau_tp_1,\dots ,\tau_tp_s\rangle_{\mathbb{Z}}$ has  rank $st$ over $\ze$.
Clearly, for $1\le i\le t$ and $1\le j\le s$ the elements $\tau_ip_j$ are   integral generators of $\Gamma_{\mathbb{Z}}$. 
Consider the  normed  space $(M\otimes_\ze \mathbb{R}, ||\cdot
||)$, in which $\Gamma_{\mathbb{Z}}$ is embedded, and endow $\Gamma_{\mathbb{Z}}$ with the induced metric.

Apply Lemma \ref{bg} to $(\Gamma_\ze,||\cdot||)$ with $n=st$. Then, there exist integral generators $\rho_1,
\dots, \rho_{st}$ of $\Gamma_{\mathbb{Z}}$ satisfying 
 \begin{equation}
 \label{stella1}
 \begin{split}
 \Big|\Big|\sum_i\alpha_i\rho_i \Big|\Big|^2&\ge  {{\czero}}(st)\sum_i|\alpha_i|^2||\rho_i||^2\\
 &\ge {{\czero}}(st)\sum_i|\alpha_i|^2\min_k||\rho_k||^2, 
 \end{split}
 \end{equation}
for all $\alpha_1,\dots , \alpha_{st}\in \ze$.

We decompose the elements $b_1, \dots, b_s \in R$ 
as $$b_i=\sum_{j=1}^t\alpha_{ij}\tau_j$$ with $\alpha_{ij}\in \ze.$
We denote 
$$\alpha=(\alpha_{11},\dots, \alpha_{1t}, \dots, \alpha_{s1},\dots, \alpha_{st})\in \ze^{st}.$$
 As usually $(\cdot)^{\rm{t}}$ indicates the transpose, we denote
$$p^{\tau}=(\tau_1p_1, \dots, \tau_tp_1, \tau_1 p_2,\dots, \tau_tp_2, \dots, \tau_1p_s, \dots, \tau_tp_s)^{\rm{t}}\in \Gamma_{\ze}^{st}$$
and $$\rho=(\rho_1,
\dots, \rho_{st})^{\rm{t}}\in \Gamma_{\ze}^{st}.$$ 
Let $P\in SL_{st}(\ze)$ be the  base change matrix such that 
$$p^{\tau}=P\rho.$$

Then $$\sum_{i}b_ip_i= \sum_{ij}\alpha_{ij}\tau_j p_i=\alpha\cdot p^{\tau}=\alpha\cdot (P\rho)=(\alpha P)\cdot  \rho.$$

Passing to the norms and using relation (\ref{stella1}) with  the coefficients $(\alpha_1,\dots , \alpha_{st})=\alpha P$, we deduce

\begin{equation*}
\Big|\Big|\sum_{i}b_ip_i\Big|\Big|^2=||(\alpha P)\cdot  \rho||^2\ge \czero(st)|\alpha P|_2^2\min_k||\rho_k||^2  ,
\end{equation*}
where $|\cdot|_2$ is the standard Euclidean norm.
On the other hand,  the triangle inequality gives
\begin{equation*}
\begin{split}
|b_i|^2_{R}& \le  \max_{k}|\tau_k|^2_R\left(\sum_{j=1}^t|\alpha_{ij}|\right)^2\\
&\le t \max_{k}|\tau_k|^2_R\sum_{j=1}^t|\alpha_{ij}|^2.
\end{split}
\end{equation*}
  We deduce
  \begin{equation*}
\frac{||\sum _ib_ip_i||^2}{\sum_i|b_i|^2_R||p_i||^2 }\ge  {\frac{\czero(st)}{t\max_j|\tau_j|^2_R}}\frac{ \min_i||\rho_i||^2}{\max_i ||p_i||^2} \frac{|\alpha P|_2^2}{|\alpha|_2^2}.
\end{equation*}
We shall still estimate $ \frac{|\alpha P|_2^2}{|\alpha|_2^2}$ independently of $\alpha$. For a linear operator $A$ and a row vector $\beta$, it holds the classical norm relation $|\beta A|_2\le H(A) |\beta|_2$. For $A=P^{-1}$ and $\beta=\alpha P$, we deduce $$ \frac{|\alpha P|_2^2}{|\alpha|_2^2}\ge \frac{1}{H(P^{-1})^2}.$$
Then
  \begin{equation*}
\frac{||\sum _ib_ip_i||^2}{\sum_i|b_i|^2_R||p_i||^2 }\ge  {\frac{\czero(st)}{t\max_j|\tau_j|^2_R}}\frac{ \min_i||\rho_i||^2}{\max_i ||p_i||^2} \frac{1}{H(P^{-1})^2}
\end{equation*}
or equivalently
\begin{equation*}
\left|\left|\sum b_i p_i\right|\right|^2\ge \cuno(p,\tau)\sum_i|b_i|^2_R||p_i||^2,
\end{equation*}
where
\begin{equation}
\label{cuno}
\cuno(p,\tau)={\frac{\czero(st)}{t\max_j|\tau_j|^2_R}}\frac{ \min_i||\rho_i||^2}{\max_i ||p_i||^2} \frac{1}{H(P^{-1})^2}.
\end{equation}

\end{proof} 
 The following non surprising proposition has some surprising
 implications; it allows us to prove Theorems \ref{altzero} and  \ref{equivalenza}.

\begin{propo}
\label{ret}

Let  $p_1, \dots, p_s$ be  linearly independent points of
$E(\overline{\qe})$ and $p=(p_1, \dots, p_s)$.  
Let $\tau$ be a set of integral generators of $\emor(E)$.
Then, there exist  positive reals $\cdue(p,\tau)$   and  $ \varepsilon_0(p,\tau)$ such that 
$$\cdue(p,\tau)\sum_i|b_i|^2||p_i||^2 \le \Big|\Big|\sum_ib_i(p_i-\xi_i)-b\upxi \Big|\Big|^2$$
 for all  $b_1,\dots,b_s,b\in \emor(E)$ with $|b|\le \max_i|b_i|$ and  for all
$\xi_1,\dots, \xi_s ,\upxi \in E(\overline{\qe})$ with $||\xi_i||, ||\upxi|| \le
 \varepsilon_0(p,\tau)$.  

In particular $p_1-\xi_1, \dots , p_s-\xi_s$ are linearly independent points of $E$.

\end{propo}

\begin{proof}
Recall that  the norm on ${\rm{End}}(E)$ is compatible
with the height norm on $E(\overline{\qe})$. Namely
$||b_ip_i||=|b_i|_{\emor(E)}||p_i||$. Thus $(\emor(E), |\cdot |)$ is a
hermitian free $\mathbb{Z}$-module of rank $1$ if $E$ has not C.M. or $2$ is $E$ has C.M. 
Furthermore,  $(E,||\cdot||)$ is a hermitian ${\rm{End}}(E)$-module.

Apply Proposition \ref{angoli} with $R=\emor (E)$, $M=E$ and $\tau=(1)$  if $\emor(E)\cong \ze$ or $\tau=(1,\tau_2)$ if $\emor(E)\cong \ze+\tau_2\ze$. We deduce that, for $b_1, \dots , b_s \in \emor(E)$
\begin{equation}
\label{nn}
 \Big|\Big|\sum b_i p_i \Big|\Big|^2\ge \cuno(p,\tau)\sum_i|b_i|^2||p_i||^2.
\end{equation}
Let $||\xi_i||, ||\upxi|| \le
\varepsilon$. Since $|b|\le \max|b_i|$ the triangle inequality implies 
\begin{equation*}
\begin{split}
 \Big|\Big|\sum_ib_i(p_i-\xi_i)-b\upxi \Big|\Big|&\ge  \Big|\Big|\sum_ib_ip_i \Big|\Big| -\varepsilon\sum_i|b_i| -
\varepsilon|b|\\
& \ge   \Big|\Big|\sum_ib_ip_i \Big|\Big|-2\varepsilon\sum_i|b_i|. 
\end{split}
\end{equation*}
Pass to the squares and do not forget that $\left(\sum^s_{i=1}|b_i|\right)^2\le s\sum^s_{i=1}|b_i|^2$. 
We deduce
\begin{equation*}
\begin{split}
 \Big|\Big|\sum_ib_i(p_i-\xi_i)-b\upxi \Big|\Big|^2& \ge
 \Big|\Big|\sum_ib_ip_i \Big|\Big|^2 -4 \varepsilon\Big|\Big|\sum_ib_ip_i \Big|\Big|\sum_i|b_i| +4 \varepsilon^2\left(\sum_i|b_i|\right)^2 \\&\ge
 \Big|\Big|\sum_ib_ip_i \Big|\Big|^2-4s\varepsilon\left(\sum_i|b_i|^2\right)\max_i||p_i||.
\end{split}
\end{equation*}
Choose
 \begin{equation}
 \label{e2}
 \varepsilon\le \varepsilon_0(p,\tau)= \frac{\cuno(p,\tau)}{8 s}\frac{\min_i||p_i||^2}{ \max_i||p_i||}.
  \end{equation}
Using relation (\ref{nn}), we deduce
\begin{equation*}
\begin{split}
\Big|\Big|\sum_ib_i(p_i-\xi_i)-b\upxi\Big|\Big|^2& \ge
\cuno(p,\tau)\sum_i|b_i|^2||p_i||^2-\frac{1}{2}\cuno(p,\tau)\left(\sum_i|b_i|^2\right)\min_i||p_i||^2\\
&\ge \frac{1}{2}\cuno(p,\tau)\sum_i|b_i|^2||p_i||^2.
\end{split}
\end{equation*}
Set for example
\begin{equation}
\label{cdue}
c_2(p,\tau)=\frac{1}{2}\cuno(p,\tau),
\end{equation} where $\cuno(p,\tau)$ is as defined in relation (\ref{cuno}).

In particular such a relation, with $b=0$, implies that only the trivial linear combination of $p_1-\xi_1, \dots , p_s-\xi_s$ is zero.

\end{proof}

At last, we  write a lemma which enables us to choose a nice maximal free set of $\Gamma_0$, the  saturated module of a sub-module $\Gamma$ of $E(\overline{\qe})$ of finite rank, as  defined in relation (\ref{gammazero}). There is nothing deep here, as we are working on finite dimensional $\co$-vector spaces. 

\begin{lem}[Quasi orthonormality]
\label{basegzero}
Let $\Gamma_0$ be the  saturated module of  $\Gamma$. Let $s$ be the rank of $\Gamma_0$.
 Then for any  real $K>0$, there exists a maximal free set $\gamma_1, \dots , \gamma_s$  of $\Gamma_0$ such that $ ||\gamma_i||\ge K$ and 
for all $b_1, \dots, b_s\in \emor(E)$
\begin{equation*}
\Big|\Big|\sum_i b_i\gamma_i\Big|\Big|^2\ge \frac{1}{9}\sum_i|b_i|^2||\gamma_i||^2.
\end{equation*}

\end{lem}
\begin{proof}
Recall that $\emor(E)$ is an order in an imaginary quadratic field $k$. Furthermore, the height norm $||\cdot||$ makes $\Gamma_0$ a hermitian $\emor(E)$-module. Let $\Gamma^{{\rm{free}}}$ be a submodule of $\Gamma_0$ isomorphic to its free part. Then $\Gamma^{{\rm{free}}}$ is a $k$ vector space of dimension $s$. Its tensor product with $\co$ over $k$ is a normed $\co$ vector space of dimension $s$, and $\Gamma^{{\rm{free}}}$ is isomorphic to $\Gamma^{{\rm{free}}} \otimes 1$. Using for instance the Gram-Schmidt orthonormalisation algorithm  in $\Gamma^{{\rm{free}}}\otimes_k \co$, we can choose an orthonormal basis $$v_i=g_i\otimes \rho_i.$$ 
So
$$\Big|\Big|\sum_ib_iv_i\Big|\Big|^2=\sum_i|b_i|^2.$$

Decompose $\rho_i=r_{i1}+\tau r_{i2}$ for $1,\tau$ integral generators of $\emor(E)$ and $r_{ij}\in \mathbb{R}$. 
Choose $\delta={(2(1+|\tau|)\max_i||g_i||)^{-1}}$
 and rationals  $q_{ij}$ such that $q_{ij}=r_{ij}+d_{ij}$ with $|d_{ij}|\le \delta$ (use the density of the rationals).
 
 Define 
 $$\gamma'_i=g_i\otimes (q_{i1}+\tau q_{i2})=(q_{i1}+\tau q_{i2})g_1\otimes 1\in\Gamma^{{\rm{free}}}\otimes 1,$$ and
 $$\delta_i=g_i\otimes(d_{i1}+\tau d_{i2}).$$
 Then
 $$v_i=\gamma'_i+\delta_i$$
 with
 $$||\delta_i||\le ||g_i|| (1+|\tau|)\delta\le \frac{1}{2}.$$
 The triangle inequality gives
   $$2\Big|\Big|\sum_ib_i\gamma'_i\Big|\Big|^2\ge \Big|\Big|\sum_ib_iv_i\Big|\Big|^2-2\Big|\Big|\sum_ib_i\delta_i\Big|\Big|^2.$$
 The orthonormality of $v_i$  and $||\delta_i||\le \frac{1}{2}$ imply that
 \begin{equation*}
 \begin{split}
 2\Big|\Big|\sum_ib_i\gamma'_i\Big|\Big|^2&\ge\sum_i|b_i|^2-2\sum_i|b_i|^2\frac{1}{4}\\
&= \frac{1}{2}\sum_i|b_i|^2.
 \end{split}
 \end{equation*}
 Finally $||\gamma'_i||\le ||v_i||+||\delta_i||\le \frac{3}{2}$, so 
 $$\Big|\Big|\sum_ib_i\gamma'_i\Big|\Big|^2\ge\frac{1}{9}\sum_i|b_i|^2||\gamma'_i||^2.$$
 It is evident that for any integer $n_0$ the same relation holds
 $$\Big|\Big|\sum_ib_in_0\gamma'_i\Big|\Big|^2\ge\frac{1}{9}\sum_i|b_i|^2||n_0 \gamma'_i||^2.$$
Let $n_0$ be an integer such that $n_0\ge {2K}$. Note that $||\gamma'_i||\ge ||v_i||-||\delta_i||\ge\frac{1}{2}$, so
$$||n_0\gamma'_i||\ge K.$$

We conclude that the maximal free set $\gamma_i=n_0\gamma_i'$ satisfies the desired conditions.

\end{proof}
Remark that we cannot directly choose  an orthonormal basis in $\Gamma^{{\rm{free}}}$, because the norm has values in $\mathbb{R}$ and not in $\qe$. 
Actually, we could prove that for any small positive real $\delta$,   there exists a maximal free set $\gamma_1,\dots,\gamma_s$ such that
$$\Big|\Big|\sum_i b_i\gamma_i\Big|\Big|^2\ge \frac{(1-\delta)^2}{(1+\delta)^2}\sum_i|b_i|^2||\gamma_i||^2.$$

We also observe that we  could use Proposition \ref{angoli} for any maximal free set  $\gamma_1,\dots,\gamma_s$ of $\Gamma_0$, and carry around the related  constant $c_1(\gamma,\tau)$.  However, we prefer absolute constants, when possible.

\section{ Gauss-reduced morphisms}

The aim of this section is to  show that  we can consider our union over Gauss-reduced algebraic subgroups, instead of the union over all algebraic subgroups. 

Let $B$ be an algebraic subgroup of $E^g$ of codimension $r$. Then $B \subset \ker \phi_B$ for a  surjective morphism $\phi_B:E^g\to E^r$.
 Conversely, we denote by $B_\phi$ the kernel of a surjection $\phi:E^g \to E^r$. Then $B_\phi$  is an algebraic subgroup of $E^g$ of codimension $r$. 
 
 The matrices  in $M_{r\times g}(\emor(E))$ of the form 
 \begin{equation*}
\phi=(aI_r|L)=
 \left(\begin{array}{cccccc}
a&\dots &0&a_{1,r+1}&\dots &a_{1,g}\\
\vdots & & \vdots & \vdots& & \vdots\\
0&\dots &a&a_{r,r+1}&\dots &a_{r,g}
\end{array}\right),
\end{equation*} 
with $ H(\phi)=|a|$ and no common factors of  the entries (up to units), will play a key role in this work. 
For $r=g$, such a morphism becomes  the identity, and $L$  shall be forgotten. These matrices have three  main advantages:
\begin{enumerate}

\item The restriction of $\phi$ to the set $E^r\times \{0\}^{g-r}$ is nothing else than the multiplication by $a$.

\item
The image of $\mathcal{O}_\varepsilon\subset E^g$ under $\phi$ is contained in the image of $\mathcal{O}_{g\varepsilon}\cap E^r\times\{0\}^{g-r}$.
Similarly, the image of $\Gamma_0^g$ under $\phi$ is contained in the image of $\Gamma_0^r\times\{0\}^{g-r}$.

\item The matrix  $\phi$ has small height compared to other matrices with same zero component of the kernel. 
\end{enumerate}

\begin{D}[Gauss-reduced Morphisms] 
\label{gr} 
We say that a  surjective morphism $\phi:E^g \to E^r$ is Gauss-reduced of rank $r$ if: 
\begin{enumerate}
\item  There exists $a\in {\rm{End}}(E)^*$ such that $aI_{r}$ is a submatrix of $\phi$, with $I_r$ the r-identity matrix,
\item $H(\phi)=|a|$, 
\item If there exists $f\in \emor(E)$ and $\phi':E^g\to E^r$ such that $\phi=f\phi'$ then $f$ is an isomorphism.
\end{enumerate}

We say that  an algebraic subgroup is Gauss-reduced if it is the
kernel of a Gauss-reduced morphism.

\end{D}

\begin{remark}
\label{a}
Note that if $\emor(E)\cong\ze$  the condition iii. in this definition simply says that the greatest common divisor of the entries of $\phi$ is $1$ and $f=\pm 1$.

Whenever we will ensure $\emor(E)\cong \ze$, we will require, in the
definition of Gauss-reduced \ref{gr} ii., the  more restrictive condition
$H(\phi)=a$,  instead of  $H(\phi)=|a|$.  Obviously $B_\phi=B_{-\phi}$, thus all
 lemmas below hold with this `up to units'-definition of Gauss-reduced.
This assumption simplifies notations.

\end{remark}

A morphisms $\phi'$, given by a reordering of the rows of a
morphism $\phi$, have the same kernel as $\phi$. Saying that $aI_r$ is a sub-matrix  of $\phi$ fixes one permutation of the rows of $\phi$.

A reordering of the columns corresponds, instead, to a permutation of the coordinates.
  Statements  will be proven for Gauss-reduced morphisms of the form $\phi=(aI|L)$. For each other reordering of the columns  the proofs are analogous. 
Since there are  finitely many permutations  of $g$ columns, the finiteness statements will follow.

The following lemma is a simple useful trick to keep in mind.
\begin{lem}
\label{ridi}

Let $\phi:E^g\to E^r$ be Gauss-reduced of rank $r$. 
\begin{enumerate}
\item
For $\xi=(\xi_1, \dots ,\xi_g)\in \mathcal{O}_\varepsilon$, there exists a point $ \xi'=(\xi'', \{0\}^{g-r})\in \mathcal{O}_{g\varepsilon}$  such that 
$$\phi(\xi)=\phi(\xi')=[a]\xi''.$$
\item
For $y=(y_1, \dots ,y_g)\in \Gamma_0^g$, there exists a point $ y'=(y'',\{0\}^{g-r})\in \Gamma_0^r\times\{0\}^{g-r}$  such that 
$$\phi(y)=\phi(y')=[a]y''.$$
\end{enumerate}
\end{lem}
\begin{proof}
Up to reordering of the columns, the morphism  $\phi$ has the form
\begin{equation*}
\phi= 
 \left(\begin{array}{cccccc}
a&\dots &0&a_{1,r+1}&\dots &a_{1,g}\\
\vdots & & \vdots & \vdots& & \vdots\\
0&\dots &a&a_{r,r+1}&\dots &a_{r,g}
\end{array}\right),
\end{equation*}
with $H(\phi)=|a|$.

i. Consider a point $\xi''\in E^r$  such that
\begin{equation*}
[a]\xi''={\phi(\xi)}.
\end{equation*}
Since  $||\xi''|| =\frac{||\phi(\xi)||}{|a|}=\max_i\frac{||\sum_j a_{ij}\xi_j||}{|a|}$ and $|a|=\max_{ij}|a_{ij}|$, we obtain  $$||\xi'' ||\le g \varepsilon.$$
Define $\xi'=(\xi'',\{0\}^{g-r})$.
 Clearly $$\phi(\xi')=[a]\xi''=\phi(\xi).$$
 
\vspace{0.4cm}

ii.  Note that, $\phi(y) \in \Gamma_0^r$. Since $\Gamma_0$ is a division group, the point 
 $y''$ such that 
\begin{equation*}
[a]y''={\phi(y)},
\end{equation*} 
belongs to $\Gamma_0^r$.
Define $y'=(y'',\{0\}^{g-r})$.
Then $$\phi(y')=[a]y''=\phi(y).$$

\end{proof}

In the following lemma, we show  that the zero components of $B_\phi$, for $\phi$ ranging over all Gauss-reduced morphisms of rank $r$, are all possible abelian subvarieties of $E^g$ of codimension $r$.  This is proven  using the classical Gauss algorithm, where the 
pivots have maximal absolute values.

\begin{lem} 
\label{tor}

Let $\psi:E^g\to E^r$ be a morphism of rank $r$. Then
\begin{enumerate}

\item
For every  $N\in \emor(E)^*$, it holds
$$B_{N\psi}\subset B_{\psi}+(E^r_{\rm{Tor}}\times \{0\}^{g-r}).$$
\item
There exists a Gauss-reduced  morphism $\phi:E^g \to E^r$ of rank $r$ such that 
$$B_{\psi}\subset B_{\phi}+(E^r_{\rm{Tor}}\times \{0\}^{g-r}).$$

\end{enumerate}

\end{lem}

\begin{proof}

i.  We show the first relation.

Let $b \in B_{N\psi}$, then $N\psi(b)=0$. So $\psi(b)=t$ with $t$ a
$N$-torsion point in $E^r$. Let $\psi_1$ be an invertible
$r$-submatrix of $\psi$. 
Up to reordering of the columns, we can suppose $\psi=(\psi_1|\psi_2)$. 
Let $t'$ be a torsion point in $E^r$ such that $\psi_1 (t')=t$.
 Then $\psi(b-(t',0))=0$. Thus  $ b \in B_{\psi}+(E^r_{\rm{Tor}}\times \{0\}^{g-r})$.

\vspace{0.4cm}

ii. We show the second relation.

The Gauss algorithm gives   an invertible integral $r$-matrix $\Delta $ such that, up to the order of the columns, $\Delta\psi$ is
of the form
\begin{equation*}
\Delta\psi= 
 \left(\begin{array}{cccccc}
a&\dots &0&a_{1,r+1}&\dots &a_{1,g}\\
\vdots & & \vdots & \vdots& & \vdots\\
0&\dots &a&a_{r,r+1}&\dots &a_{r,g}
\end{array}\right),
\end{equation*}
with $H(\Delta\psi)=|a|$ (potentially there are common factors of the entries).

Let $b\in B_\psi$, then $\psi(b)=0$. Hence $\Delta\psi(x)=0$. It follows
$$
B_\psi \subset B_{\Delta\psi}.$$
Let $N\in \emor(E)^*$ such that  $N|\Delta\psi$ and such that   if $f|(\Delta\psi/N)$ then $f$ is a unit (if $\emor(E)\cong \ze$, then $N$ is simply  the greatest common divisor of the entries of $\Delta\psi$).
We define $$\phi=\Delta\psi/N.$$ Clearly $\phi$ is Gauss-reduced and $B_\psi \subset B_{\Delta\psi}=B_{N\phi}$. By part i. of this lemma applied to $N\phi$, we conclude  $$B_\psi \subset  B_{\phi}+(E^r_{\rm{Tor}}\times \{0\}^{g-r}).$$






\end{proof}
Note that, in the previous lemma, a reordering of the columns of $\psi$ or $\phi$   induces the same reordering of the coordinates of $E^r_{{\rm{Tor}}}\times\{0\}^{g-r}$.

Taking intersections with the algebraic points of our curve, the previous lemma part ii.  translates  immediately as
\begin{lem}
\label{tor1}
Let  $C\subset E^g$ be an algebraic curve (transverse or not).
For any real $\varepsilon\ge0$
$$S_r(C,(\Gamma_0^g)_\varepsilon)=\bigcup_{\substack{\phi\,\,\,\, {\rm{Gauss-reduced}}\\
{\rm{rk}}(\phi)= r}} C(\overline{\qe})\cap (B_\phi+(\Gamma_0^g)_{\varepsilon} ).$$

\end{lem}
\begin{proof}
By definition
$$S_r(C,(\Gamma_0^g)_\varepsilon)\supseteq \bigcup_{\substack{\phi\,\,\,\, {\rm{Gauss-reduced}}\\
{\rm{rk}}(\phi)= r}}C(\overline{\qe})\cap  (B_\psi+(\Gamma_0^g)_{\varepsilon} ).$$
On the other hand, 
by the previous Lemma \ref{tor} ii, we see that
$$C(\overline{\qe})\cap  (B_\psi+(\Gamma_0^g)_{\varepsilon} ) \subset  C(\overline{\qe})\cap (B_\phi+(E^r_{{\rm{Tor}}}\times \{0\}^{g-r})+(\Gamma_0^g)_{\varepsilon} ),$$
with $\phi$ Gauss-reduced of rank $r$.
Moreover $(E^r_{{\rm{Tor}}}\times \{0\}^{g-r})\subset \mathcal{O}_{\varepsilon}\subset (\Gamma_0^g)_{\varepsilon}$.
\end{proof}

\section{Relation between transverse and  weak-transverse curves}

We discuss here how  we can associate to  a couple $(C,\Gamma)$, with  $C$ a transverse curve and $\Gamma$ a subgroup of finite rank, 
a weak-transverse curve $C'$  and vice versa.  There are properties which are easier for $C$ and others for $C'$. Using this association, we will try to gain advantages from both situations. 


\subsection{From transverse to weak-transverse}
Let $C$ be   transverse  in $E^g$.  If $\Gamma$ has rank $0$, we set $C'=C$.  If ${\rm{rk}}\,\,\Gamma\ge 1$, consider   the saturated module $\Gamma_0$ of rank $s$ associated to $\Gamma$, as defined in   relation (\ref{gammazero}).  Let $\gamma_1, \dots , \gamma_s$ be a maximal free set of $\Gamma_0 $. 
 We denote the associated point of $E^s$  by $$\gamma=(\gamma_1, \dots , \gamma_s).$$
 We define
 $$C'=C\times \gamma.$$
  Since  $C$
  is transverse and the $\gamma_i$ are $\emor(E)$-linearly independent, the curve $C'$ is weak-transverse.  More precisely, suppose on the contrary that $C'$
  would be  contained in an algebraic
subgroup $B_\phi$ of codimension $1$, with $\phi=(a_1,\dots ,a_{g+s})$. Define
$y_1$ to be a point in $E$ such that  $a_1y_1=\sum_{i=g+1}^{g+s}a_i\gamma_{i-g}$
and define $y=(y_1,0,\dots,0)\in E^g$. Then $C\subset B_{\phi_1}+y$ with $\phi_1=(a_1,\dots,a_g)$, contradicting that $C$ is transverse.

\subsection{From weak-transverse to transverse}
Let $C'$ be weak-transverse in $E^{n}$.  If $C'$ is transverse, we set $C=C'$ and $\Gamma=0$. Suppose that $C'$ is not transverse.
Let $H_0$ be the abelian
subvariety of smallest dimension $g$ such that $C'\subset H_0+p$ for
$p\in H_0^\perp(\overline\qe)$ and $H_0^\perp$  the orthogonal complement  of $H_0$ with respect to the canonical polarization.
 
 Then  $E^n$ is isogenous to $ H_0\times H_0^{\perp}$. Furthermore $H_0$
 is isogenous to $E^g$ and $H_0^{\perp}$ is isogenous to $E^s$ where $s=n-g$.  Let
 $j_0$, $j_1$ and $j_2$ be such  isogenies. 
 We fix the isogeny $$j=(j_1\times j_2)\circ j_0:E^n \to H_0\times H_0^{\perp}\to E^g\times E^s,$$
which sends $H_0$ to $E^g\times 0$ and $H_0^{\perp}$ to $0\times E^{s}$.

Then $$j(C')\subset (E^g\times 0) +j(p),$$ with $j(p)=(0,
\dots ,0, p_1,\dots,p_s)$.

We consider the natural projection on the first $g$ coordinates
\begin{equation*}
\begin{split}
 \pi:&E^g\times E^s \to E^g\\
 &j(C') \to \pi(j(C')).
 \end{split}
 \end{equation*}
 We define
$$C=\pi(j(C'))\,\,\,\,\,{\rm{and}}\,\,\,\, \Gamma=\langle p_1,\dots,p_s\rangle^g.$$
 Since $H_0$ has minimal dimension, the curve $C$ is transverse in $E^g$.
 
 Note that $$j(C')=C\times (p_1, \dots, p_s).$$
 In addition  $j(C')$ is weak-transverse, because $C'$  is. Therefore,
  $\langle p_1, \dots , p_s\rangle$ has rank $s$; indeed if $\sum_{i=1}^sa_ip_i=0$,
  then $j(C')\subset B_{\phi}$ for $\phi=(\{0\}^{g},a_1,\dots ,a_s)$.

\subsection{ Weak-transverse up to an isogeny}
\label{trwtr}
Statements on boundedness of heights or finiteness of sets  are  invariant under an
isogeny of the ambient variety. Namely, given an isogeny $j$ of $E^g$,  Theorems  \ref{altzero} and \ref{tutto} hold for a curve if and only if they hold for its image via $j$. Thus,
the previous discussion shows that without loss of generality, we can assume that a weak-transverse curve $C'$ in $E^{n}$ is of the form
$$C'=C\times p$$ with 
\begin{enumerate}
\item $C$ transverse in $E^g$, 
\item    $p=(p_1, \dots ,p_s)$  a point in $E^s$ such that the module  $\langle p_1, \dots
  ,p_s\rangle$ has rank $s$,
  \item $n=g+s$.
\end{enumerate}
 This simplifies the setting for weak-transverse curves.

\subsection{Implying  the Mordell-Lang plus Bogomolov Theorem for curves}
\label{mlb}
Note that $$S_g(C,\mathcal{O}_\varepsilon)=C\cap \mathcal{O}_\varepsilon$$ and 
$$S_g(C(\Gamma_0^g)_\varepsilon)=C\cap (\Gamma_0^g)_\varepsilon.$$Moreover 
$S_2(C,\cdot)\supset S_g(C,\cdot)$. This immediately  shows that Conjecture \ref{geni} implies the Bogomolov Theorem for weak-transverse curves and the Mordell-Lang plus Bogomolov Theorem for transverse curves. We want to show that Conjecture \ref{geni} implies these theorems for all curves of genus $\ge 2$. 

In $E^g$ a curve of genus $2$ is a translate of an elliptic curve isogenous to $E$.
If $C$ is not transverse, then $C\subsetneq H_0+p$ with $H_0$ an algebraic subgroup of minimal dimension satisfying  such inclusion.  Let $\pi:E^g \to E^g/H^\perp_0$  be the natural projection and  let $\psi:E^g/H^\perp_0 \to E^k$ be an isogeny. Then $||\psi \pi(x)|| \ll ||x||$. 
In $E^k$, consider the transverse curve $C' =\psi \pi(C-p) $ and $\Gamma'=\psi\pi \langle \Gamma, \Gamma_p\rangle$.  Note that  $\psi \pi({\rm{Tor}}_{E^g} )\subset {\rm{Tor}}_{E^k}$. Then 
$$S_g(C, (\Gamma_0^g)_\varepsilon)\subset  \pi^{-1}_{|C} S_k(C',  ({\Gamma'}_0^g)_{\varepsilon'}).$$
The map $\pi^{-1}_{|C}$ has  finite fiber. 
Applying Conjecture \ref{geni} to $C'\subset E^k$ we deduce that $S_g(C, (\Gamma_0^g)_\varepsilon)$ is finite.

Note that such a proof works only for $S_g(C, \cdot)$.
 Indeed the projection $\psi\pi(B)\subset E^k$
  of an algebraic subgroup $B$ of $E^g$ 
  of codimension $r$, might not be of codimension
   $r$ in $E^k$. Eventually, it could  be the entire $E^k$.

\section{Quasi-special Morphisms}
As Gauss-reduced morphisms play a key role for transverse curves,
quasi-special morphisms play a key role for weak-transverse curves. In
particular, for small $\varepsilon$,  quasi-special morphisms will be
enough to cover the whole $S_r(C\times p,\mathcal{O}_\varepsilon)$ -
see Lemma \ref{ind} below.

 Let us give a flavor for quasi-special. Suppose that $C\times p$ is weak-transverse in $E^{g+s}$ with $C$  transverse in $E^g$. A point of $C\times p$ is of the form $(x,p)$. The last $s$-coordinates are constant and just the $x$ varies. This two parts must be treated differently. Saying that a morphism $\tilde\psi=(\psi|\psi')$ is quasi-special ensures that the rank of $\psi$ is maximal (note that $\psi$ acts on $x$).  In particular, this allows us to apply the Gauss algorithm on the first $g$ columns of $\tilde\phi$.
\begin{D}[Quasi-special]
A surjective morphism  $\tilde{\phi}:E^{g+s} \to E^r$ is quasi-special if
 there exist  $N\in \zono^*$,  morphisms  $\phi:E^g\to E^r$ and $\phi':E^s \to E^r$ such that 
\begin{enumerate}
\item$\tilde{\phi}=(N\phi|\phi')$,
\item  $\phi=(aI_r|L)$ is Gauss-reduced of rank $r$,
\item   If there exists $f\in \emor(E)$ and $\tilde\phi':E^{g+s}\to E^r$ such that $\tilde\phi=f\tilde\phi'$ then $f$ is an isomorphism.
\end{enumerate}
\end{D}
Note that we do not require that $\tilde\phi$ is Gauss-reduced, the
fact is that  $H(\phi')$ might not be controlled by $NH(\phi)$. This extra condition will define special morphisms (see Definition \ref{spec}).

\begin{lem}
\label{ind}
Let $C\times p$ be weak-transverse in $E^{g+s}$ with  $C$ transverse in $E^g$. 
Then, there exists $\varepsilon>0$ such that 
$$S_r(C\times p,\mathcal{O}_\varepsilon)\subset
\bigcup_{\substack{\tilde\phi\,\,\,  {\rm{quasi-special}}\\{\rm{rk}}\tilde\phi=r}} (C(\overline{\qe})\times p)\cap (B_{\tilde\phi}+\mathcal{O}_{\varepsilon} ).$$
We can choose $\varepsilon\le  \ezero$, where $\ezero$ is as in
Proposition \ref{ret}.

\end{lem}
\begin{proof}
 Let $(x,p)\in S_r(C\times p ,\mathcal{O}_\varepsilon)$. Then $(x,p) \in (C(\overline{\qe})\times p)\cap (B_{\tilde\psi}+\mathcal{O}_{\varepsilon} )$ for  a morphism $\tilde\psi=(\psi|\psi'):E^{g+s}\to E^r$ of rank $r$. In  other words, there exists a point $(\xi,\xi')\in \mathcal{O}_\varepsilon$ such that $$\tilde\psi((x,p)+(\xi,\xi'))=0.$$

 First, we show that the rank of $\psi$ is $r$.
 Suppose, on the contrary, that the rank of $\psi$ would be less than $r$. Then  a linear combination  of the rows of $\psi$ is trivial, namely $$(\lambda_1,\dots ,\lambda_r )\psi=0.$$
 Since $\psi(x+\xi)+\psi'(p+\xi')=0$, the same linear combination of the  $r$ coordinates of $\psi'(p+\xi')$ is trivial, namely $$(\lambda_1,\dots ,\lambda_r )\psi'(p+\xi')=0.$$
 Apply Proposition \ref{ret} with $ (b_1,\dots,b_s)=(\lambda_1,\dots ,\lambda_r )\psi'$, $(\xi_1,\dots,\xi_s)=-\xi'$, $\upxi=0$ and $b=0$. This implies  that, if $\varepsilon\le  \ezero$, then the points $p_1+\xi'_1,\dots ,p_s+\xi'_s $ are linearly independent. 
 It follows that $$(\lambda_1,\dots ,\lambda_r)\psi'=0.$$ Hence, the rank of $\tilde\psi$ would be  less than $r$, contradicting the fact that the rank of $\tilde\psi$ is $r$.

Since the rank of $\psi $ is $r$, we can apply the Gauss algorithm using pivots in  $\psi$ of maximal absolute values in $\psi$ (clearly we cannot require that they have maximal absolute values in $\tilde\psi$).
Let $\Delta$ be an invertible matrix, given by the  Gauss algorithm, such that $\Delta \tilde\psi=(\phi_1|\phi_2)$ with $fI_r$ a submatrix of $\phi_1$.

We shall still get rid of possible common factors.
Let $N_1, n_1 \in \emor(E)^*$ such that  $N_1|\phi_1$ and $n_1|\Delta\tilde\psi$. Further suppose  that, if $f|(\phi_1/N_1)$ or
$f|(\Delta\tilde\psi/n_1)$ 
then $f$ is a unit of $\emor(E)$ (if $\emor(E)\cong \ze$, then $N_1$ is the greatest
common divisor of the entries of $\phi_1$ and $n_1$  the greatest
common divisor of the entries of $\Delta\tilde\psi$ ).
Then   $$\Delta\tilde\psi=n_1(N\phi|\phi')$$
with $N=N_1/n_1$, $\phi=\phi_1/N_1$ and $\phi'=\phi_2/n_1$.
We define
$$\tilde\phi=(N\phi|\phi').$$
Clearly $\tilde\phi$ is quasi-special.
In addition $$B_{\tilde\psi}\subset B_{\Delta\tilde\psi}= B_{n_1\tilde\phi}.$$
By Lemma \ref{tor} i. (with $\psi=\tilde\phi$ and $N=n_1$) we deduce that 
$$B_{\tilde\psi}\subset B_{\tilde\phi}+E^r_{\rm{Tor}}\times\{0\}^{g+s-r}.$$
Since $(x,p)\in B_{\tilde\psi}+\mathcal{O}_\varepsilon$, then
$(x,p)\in B_{\tilde\phi}+\mathcal{O}_\varepsilon$ with $\tilde\phi$ quasi-special. 
\end{proof}

\section{Estimates for the Height: the Proof of Theorem \ref{altzero}}
\label{bound}

As it has been already pointed out, Theorem \ref{altzero} part ii. is proven in \cite{RV}  Theorem
1.5.
In this section, we adapt the proof of  \cite{RV}  Theorem 1.5 
to Theorem \ref{altzero} part i.

In view of  section \ref{trwtr},  we can assume, without loss of generality, that  a weak-transverse curve $C'$ in $E^{n}$ has   the form
$$C'=C\times p$$ with 
\begin{enumerate}
\item $C$ transverse in $E^g$, 
\item   $p=(p_1, \dots ,p_s)$ a point in $E^s$ such that the  module  $\langle p_1, \dots
  ,p_s\rangle$ has rank $s$,
  \item $n=g+s$.
\end{enumerate}


 \begin{D}
 
 Let $p$ be a point in $E^s$ and $\varepsilon$ a non negative real.
    We define  $G_p^\varepsilon$ as the set of points  $\g \in E^2$ for which there exist a matrix $A\in M_{2, s} (\zono)$, an element $a\in \zono$ with  $0<|a|\le H(A)$, points $\xi \in E^s$ and $\upxi \in E^2$ of norm at most $\varepsilon$ such that 
    $$[a]\g={A}({p}+\xi) +[a]\upxi.$$
We identify $G_p^{\varepsilon}$ with the subset $G_p^{\varepsilon}\times\{0\}^{g-2}$ of $E^g$.  
\end{D}
Now we embed $S_2(C\times {p},\mathcal{O}_{\varepsilon})$ in two sets related to the transverse curve $C$. We then use the Vojta inequality on these new sets.
\begin{lem}
\label{spezzo}
The natural projection  on the first $g$ coordinates 
\begin{equation*}
\begin{split}
E^{g}\times E^s &\to E^g\\
(x,y) &\to x
\end{split}
\end{equation*}
 defines an injection
\begin{equation*}
 S_2(C\times {p},\mathcal{O}_{\varepsilon/2gs})\hookrightarrow S_2(C,(\Gamma^g_p)_\varepsilon)\cup \bigcup_{\substack{\phi:E^g\to E^2 \\{\rm{Gauss-reduced}}}}C(\overline{\qe})\cap B_\phi+G_p^\varepsilon .
 \end{equation*}
\end{lem}

\begin{proof}

Let $(x,{p}) \in S_2(C\times {p},\mathcal{O}_{\varepsilon/2gs})$. By Lemma \ref{ind}, $(x,p)\in B_{\tilde\phi}+\mathcal{O}_{\varepsilon/2gs}$,  with ${\tilde\phi}=({N\phi}|{\phi'}):E^{g+s} \to E^2$ quasi-special of rank $2$. Hence 
 $${\tilde\phi}((x,{p})+(\xi,\xi'))=0,$$
 for $(\xi,\xi')\in \mathcal{O}_{\varepsilon/2gs}$.   
 We can write the equality as
 $$N\phi(x)+N\phi(\xi)+\phi'(p+\xi')=0.$$  
 By definition of quasi-special $\phi$ is Gauss-reduced, so  $${\phi}=(aI_2|L).$$
  By Lemma \ref{ridi} i applied to $\phi$ and $\xi$, we can assume $$\xi=(\xi_1, \xi_2, 0,\dots 0)\in \mathcal{O}_{\varepsilon/2s}.$$

Suppose first  that $NH({\phi})\ge H({\tilde\phi})$. 
Let $\upxi$ be a point in $E^2\times \{0\}^{g-2}$ such that
 $$N[a]\upxi=({\phi'}(\xi'), 0\dots,0).$$ Then $$||\upxi||=\frac{{||\phi'(\xi')||}}{NH({\phi})}\le \frac{\varepsilon}{2}.$$ 
 Let $y$ be a point in $E^2\times \{0\}^{g-2}$ such that 
  $$N[a]y=( {{\phi'}}({p}),0,\dots ,0).$$
  Since $\Gamma_p$ is a division group, $y \in \Gamma_p^2\times \{0\}^{g-2}$.
  Then $$N{\phi}(x+\xi+\upxi+y)=0$$ with $y +\xi+\upxi\in \Gamma_p^g+\mathcal{O}_\varepsilon$.
   So 
  $$x\in S_2(C,(\Gamma_p^g)_{\varepsilon}).$$ 

Suppose now  that $NH({\phi})< H({\tilde\phi})$ or equivalently  $NH({\phi})< H({\phi'})$. 
Let $\g'$ be a point in $E^2$ such that 
$$N[a]\g' ={{\phi'}}({p}+\xi')+N[a] {\xi_1 \choose \xi_2},$$  and $\g=(\g',\{0\}^{g-r})$. Then $\g\in G^{\varepsilon}_p$.
Moreover 
\begin{equation*}
\begin{split}
N{\phi}(x+\g)&=N\phi(x)+N\phi(\g)\\
&=N\phi(x)+N[a]\g'\\
&=N\phi(x)+{{\phi'}}({p}+\xi')+N[a] {\xi_1 \choose \xi_2}\\
&=N\phi(x)+N\phi(\xi)+{{\phi'}}({p}+\xi')\\
&=\tilde\phi((x,p)+(\xi,\xi'))=0.
\end{split}
\end{equation*}
Thus $$x\in B_{N\phi} +G_p^\varepsilon,$$
and by Lemma \ref{tor} i.
$$ x \in B_\phi+(E^2_{\rm{Tor}}\times \{0\}^{g-2})+G_p^\varepsilon.$$
 Note that $G_p^\varepsilon+(E^2_{\rm{Tor}}\times \{0\}^{g-2})\subset G_p^\varepsilon$. Hence,  $$x \in C(\overline{\qe})\cap B_\phi+G_p^\varepsilon.$$

\end{proof}
  
  \begin{lem} [Equivalent of \cite{RV} Lemma 3.2]
  \label{3.2}
   For $\phi:E^g \to E^2$ Gauss-reduced of rank $2$, we have the following inclusion of sets 
\begin{equation*}
(B_\phi+G_p^\varepsilon)  {\subset}
\{ P+ \g \,\,:\,\,P\in B_\phi,\,\,\g\in G_p^\varepsilon  \,\,\,{\rm{and}}\,\,\,\max(||\g||,||P||)\le 2g||P+\g||\}.
\end{equation*}
 
 \end{lem}
 \begin{proof}
Let $x\in(B_\phi+G_p^\varepsilon)$ with $\phi=(aI_r|L)$ Gauss reduced of rank 2. Then $x=P+\g$ with  $P\in B_\phi$ and $\g \in G_p^\varepsilon$ and  $\phi(x-\g)=0$. By definition $G_p^\varepsilon\subset E^2\times \{0\}^{g-2}$, so $\phi(\g)=   [a]\g$. Then $$||\g||=\frac{||\phi(\g)||}{H(\phi)}=\frac{||\phi(x)||}{H(\phi)}\le g||x||.$$ So $$||P||=||x-\g||\le (g+1)||x||=(g+1)||P+\g||.$$



\end{proof}

Note that, \cite{RV} Lemma 3.3 part (1) is a statement on the morphism, therefore it holds with no need of any remarks.

\begin{lem}[Equivalent of \cite{RV} Lemma 3.3 part (2)]
\label{3.3}
There exists an effective $\varepsilon_2>0$ such that, for all $\varepsilon\le \varepsilon_2$,   any sequence of elements in $G_p^\varepsilon$ admits  a sub-sequence in which every two elements $\g$, $\g'$ satisfy 
$$\left|\left|\frac{\g}{||\g||}-\frac{\g'}{||\g'||}\right|\right|\le \frac{1}{16 g c_1},$$
where $c_1$ depends on $C$ and  is as  defined in \cite{RV} Proposition 2.1.
\end{lem}
\begin{proof}
We decompose two elements $\g$ and $\g'$ in a given sequence of elements  of $G_p^\varepsilon$ as follows 
$$[a]\g=A({p}+\xi)+[a]\upxi$$ and $$[a']\g'=A'({p}+\xi')+a'\upxi'$$ with $A,A'\in {M}_{2,s}(\zono)$ and $0<|a|\le H(A)$, $0<|a'|\le H(A')$. Let us define $y$ and $y'$ such that  $$[a]y=A({p})$$ and $$[a']y'=A'({p}).$$
Since the sphere of radius $1$ is compact in $$\left(\langle p_1,\dots ,p_s\rangle\times \langle p_1,\dots ,p_s\rangle\right)\otimes \mathbb{R},$$ we can extract a sub-sequence such that for every two elements $y$ and $y'$ it holds 
$$\left|\left|\frac{y}{||y||}-\frac{y'}{||y'||}\right|\right|\le \frac{1}{48gc_1}.$$
Note that
$$\left|\left|\frac{\g}{||\g||}-\frac{y}{||\g||}
  \right|\right|=\left|\left|\frac{A(\xi)+[a]\upxi}{A({p}+\xi)+[a]\upxi}\right|\right|$$
and
$$\left|\left|\frac{y}{||\g||}-\frac{y}{||y||}\right|\right|=\left|\frac{
    ||A(p)||-||A(p+\xi)+[a]\upxi||}{||A(p+\xi)+[a]\upxi||}\right| \le
\left|\left|\frac{A(\xi)+[a]\upxi}{A({p}+\xi)+[a]\upxi}\right|\right|$$ 
(same relations with  $'$).
 We deduce
\begin{equation*}
\begin{split}
\left|\left|\frac{\g}{||\g||}-\frac{\g'}{||\g'||}  \right|\right| \le 
 \left|\left|\frac{y}{||y||}-\frac{y'}{||y'||}\right|\right|&+ \left|\left|\frac{y}{||\g||}-\frac{y}{||y||}\right|\right|+ \left|\left|\frac{y'}{||\g'||}-\frac{y'}{||y'||}\right|\right|\\&+ \left|\left|\frac{\g}{||\g||}-\frac{y}{||\g||}
  \right|\right|+\left|\left|\frac{\g'}{||\g'||}-\frac{y'}{||\g'||}
  \right|\right|\\
 \le \left|\left|\frac{y}{||y||}-\frac{y'}{||y'||}\right|\right|&+2\left|\left|\frac{A(\xi)+[a]\upxi}{A({p}+\xi)+[a]\upxi}\right|\right|+2\left|\left|\frac{A'(\xi')+[a']\upxi'}{A'({p}+\xi')+[a']\upxi'}\right|\right|.
\end{split}
\end{equation*}

 Choose \begin{equation}
  \label{e3}\varepsilon \le \varepsilon_2= \min(\ezero,\etre),  \end{equation} where
  $\ezero$ is defined in relation (\ref{e2}), $\cdue(p,\tau)$ is defined in relation (\ref{cdue}) and 
  $$\etre=\frac{c_2(p,\tau)^{\frac{1}{2}}{\min}||{p}_i||}{96(s+1)c_1} .$$
Note that $||A({p}+\xi)+[a]\upxi||=||A_{k}({p}+\xi)+a\upxi_k||$ for $k=1$ or $2$ and $A={A_1\choose A_2}$. Proposition \ref{ret} applied with $b_1,\dots,b_s=A_k$, $\xi=-\xi$, $\upxi=-\upxi_k$ and $b=a$, implies
$$||A({p}+\xi)+[a]\upxi||\ge  H(A)c_2(p,\tau)^{\frac{1}{2}}{\min}||{p}_i||$$ (same relation with $'$).

 It follows

\begin{equation*}
\begin{split}
\left|\left|\frac{\g}{||\g||}-\frac{\g'}{||\g'||}\right|\right|\le \frac{1}{48 g c_1} &+\varepsilon\frac{2H(A)(s+1)}{H(A)c_2(p,\tau)^{\frac{1}{2}}\min||{p}_i||}\\&+\varepsilon\frac{2H(A')(s+1)}{H(A')c_2(p,\tau)^{\frac{1}{2}}\min||{p}_i||}\\
\le  \frac{1}{48 g c_1} &+ \frac{1}{48 g c_1} +\frac{1}{48 g c_1}, 
\end{split}
\end{equation*}
where in the last inequality we use $\varepsilon\le \etre $.
\end{proof}

We are ready to conclude.

\begin{proof}[Proof of Theorem \ref{altzero} i.]
In view of  Lemma \ref{spezzo}, we shall  prove that  there exists $\varepsilon>0$ such that 
$S_2(C,(\Gamma_p)_\varepsilon)$ and
 $ \bigcup_{\substack{\phi:E^g\to E^2 \\{\rm{Gauss-reduced}}}}
 C(\overline{\qe})\cap B_\phi+G_p^\varepsilon$ have bounded height.

By Theorem \ref{altzero} ii., there exists $\varepsilon_1>0$ such that for $\varepsilon\le \varepsilon_1$,
the first set has bounded height.

It remains to show that there exists $\varepsilon_2>0$ such that, for $\varepsilon\le \varepsilon_2$, the set 
$$ \bigcup_{\substack{\phi:E^g\to E^2 \\ {\rm{Gauss-reduced}}}}C(\overline{\qe})\cap B_\phi+G_p^\varepsilon$$ has bounded height.
The proof follows, step by step, the proof of 
 \cite{RV}  Theorem 1.5.
 In view of Lemma \ref{3.2} and \ref{3.3}, all conditions for the proof of \cite{RV} Theorem 1.5 are satisfied.
 The proof is then exactly equal to the proof \cite{RV} page 1927-1928. 
  \end{proof}

\begin{remark}

In   \cite{RV}  Theorem 1.5, we show  that for $\varepsilon_1=
\frac{1}{2^gc_1}$ the set $S_2(C,\Gamma_{\varepsilon_1})$ has bounded height.
 The
 constant  $c_1$ depends on the invariants of the curve $C$. This
 constant is defined in \cite{RV} Proposition 2.1 and it is effective. On the other hand, the height of $S_2(C,\Gamma_{\varepsilon_1})$ is bounded by a constant  which is not known to be effective, unless $\Gamma$ has rank $0$.

 For $C\times p $, we have shown that for
  $\varepsilon'_2=\frac{\min(1,c_2(p,\tau))\min||{p}_i||^2}{2^8g(s+1)^2\max||p_i||c_1}$  the set  $S_2(C\times p,\mathcal{O}_{\varepsilon'_2})$ has bounded height - see relation (\ref{e3}) and Lemma \ref{spezzo}. 
 As in the previous case, the height of $S_2(C\times p,\mathcal{O}_{\varepsilon'_2})$ is bounded by a constant which, in general, is not known to be effective.

\end{remark}

\section{Recap}
We would like to  recall and  fix the notations for the rest of the
article.

For simplicity, we assume that $\emor(E)\cong \ze$. In this case the saturated module of a group coincides with its division group. According to
Remark \ref{a}, we use $H(\phi)=a$ in the definition of a Gauss-reduced morphism and $N\in \mathbb{N}^*$ in the definition of quasi-special.

\begin{itemize}
\item Let $E$ be an elliptic curve without C.M. over $\overline{\qe}$.
\item Let $C$ be a transverse curve in $E^g$ over $\overline{\qe}$. 
\item Let 
\begin{equation*}
\phi=   \left(\begin{array}{c}
\phi_1\\
\vdots\\
\phi_r
\end{array}\right)=
 \left(\begin{array}{cccccc}
a&\dots &0&L_1\\
\vdots & &\vdots & \vdots \\
0&\dots &a&L_r
\end{array}\right)
\end{equation*}
 be a Gauss-reduced morphism of rank $1 \le r\le g$, with $L_i\in\mathbb{Z}^{g-r}$ and $H(\phi)=a$.
\vspace{0.5 cm}
\item Let $\Gamma$ be  a subgroup of finite rank of  $E^g(\overline\qe)$.
\item Let $\Gamma_0$ be the  division group of
  $\Gamma$ and $s$ its rank (the definition is given in relation (\ref{gammazero})). 
\item Choose $\euno>0$ so that $S_2(C, (\Gamma^g_0)_{\euno})$ has
  bounded height (the definition is consistent in view of Theorem
  \ref{altzero} ii.).
\item  Let $\kuno$ be the norm of   $S_2(C, (\Gamma^g_0)_{\euno})$.
\vspace{0.5 cm}
\item Let $\gamma=(\gamma_1,\dots , \gamma_s)$ be a  point of $E^s(\overline\qe)$ such that $\gamma_1,\dots , \gamma_s$ is a maximal free set of $\Gamma_0$ satisfying the conditions of Lemma \ref{basegzero} with $K=3g\kuno$. Namely, for all integers $b_i$
 \begin{equation}
 \label{con1}
 \frac{1}{9}\sum_i|b_i|^2||\gamma_i||^2\le \Big|\Big|\sum_ib_i\gamma_i\Big|\Big|^2
 \end{equation}
 and
\begin{equation}
\label{bignor}
\min_i||\gamma_i|| \ge {3g\kuno}.
\end{equation}

\item Let  $C\times \gamma$ be the associated weak-transverse curve  in $E^{g+s}$.
\item Let $\tilde\phi=(N\phi|\phi'):E^{g+s}\to E^r$ be a quasi-special morphism with $N\in \mathbb{N}^*$.

\item Choose  $\edue>0$  so that $S_2(C\times \gamma,
  \mathcal{O}_\edue)$ has bounded height (the definition is consistent in view of Theorem
  \ref{altzero} i.).
\item  Let $\kdue$ be the norm of   $S_2(C\times \gamma,
  \mathcal{O}_\edue)$.
\vspace{0.5 cm}
\item Let $p=(p_1,\dots ,p_s)\in E^s$ be a point such that  the rank of $\langle p_1,\dots ,p_s\rangle$ is $s$.
\item Let  $\Gamma_p$  be the division group  of $\langle p_1,\dots ,p_s\rangle$ (in short the division group of $p$). 

\item Let $\cunop$ and $\ezerop$ be the constants  $(\cdue(p,\tau))^{\frac{1}{2}}$ and
  $\ezero$ defined in Proposition
  \ref{ret} for the point $p$ and $\tau=1$ (please note the square root in $\cunop$).

\item Let  $C\times p$ be the associated weak-transverse curve  in $E^{g+s}$.
\item Choose $\eduep>0$ so that $S_2(C\times p, \mathcal{O}_\eduep)$
  has bounded height (the definition is consistent in view of Theorem
  \ref{altzero} i.).
\item Let  $\kduep$ be the norm of $S_2(C\times p, \mathcal{O}_\eduep)$. 

\end{itemize}

\section{Equivalence of the strong Statements: The Proof of Theorem \ref{equivalenza}}
\label{equi}
The following theorem implies Theorem \ref{equivalenza}  immediately;
in addition it gives explicit inclusions.
Once more, we would like to emphasise that we need to assume that
$S_r(C\times p, \mathcal{O}_\varepsilon)$ has bounded height in order
to embed it in a set of the type $S_r(C, \Gamma_{\varepsilon'})$. Therefore
we assume $r\ge2$ and $\varepsilon\le \eduep$ in part ii.

\begin{thm}
\label{equi34}
Let $\varepsilon\ge0$. Then,

\begin{enumerate}
\item
 The map $x \to (x,\gamma)$ defines an injection
$$S_r(C,\Gamma_{\varepsilon})\hookrightarrow  S_r(C\times \gamma,\mathcal{O}_{\varepsilon}).$$ 
Recall that $\gamma$ is a maximal free set of $\Gamma_0$.
\item
 For $2\le r $ and  $ \varepsilon \le \min(\ezerop,\eduep)$, the map $(x,p) \to x$ defines an injection
$$S_r(C\times p,\mathcal{O}_{\varepsilon})\hookrightarrow
S_r(C,(\Gamma^g_p)_{\varepsilon\ktre}),$$ where
$\ktre=(g+s)\max\left(
  1,\frac{g(\kduep+\varepsilon)}{{\cunop}\min_i||p_i||}\right)$. 
  Recall that $\Gamma_p$ is the division group of $p$.
  \end{enumerate}
\end{thm}

\begin{proof}
 i.  Let $x \in S_r(C,\Gamma_{\varepsilon})$. Then, there exists a surjective
 $\phi:E^g \to E^r$,  points $y\in \Gamma$ 
and $\xi \in \mathcal{O}_{\varepsilon}$ such that 
$$\phi(x+y+\xi)=0.$$
Since $\gamma=(\gamma_1, \dots ,\gamma_s)$ is a maximal free set of $\Gamma_0$,
 there exists  a positive integer $N$ and a matrix
$G\in M_{r,s}(\mathbb{Z})$ such that
$$[N]y=G\gamma.$$
We define $$\tilde\phi=(N\phi|\phi G).$$
Then
\begin{equation*}
\begin{split}
\tilde\phi((x,\gamma)+(\xi,0))&=N\phi(x+\xi)+\phi G(\gamma)\\
&=N\phi(x+\xi+y)=0.
\end{split}
\end{equation*}
So $$(x,\gamma) \in S_r(C\times \gamma, \mathcal{O}_\varepsilon).$$

\vspace{0.4cm}
 ii. 
Let $(x,p) \in S_r(C\times p,\mathcal{O}_{\varepsilon})$.
Thanks to Lemma \ref{ind}, the assumption $\varepsilon \le \ezerop$ implies
$$(x,p) \in (B_{\tilde\phi}+\mathcal{O}_{\varepsilon})$$ with
$\tilde\phi=(N\phi|\phi')$ quasi-special. Hence  $$\tilde\phi((x,p)+(\xi,\xi'))=0,$$ for $(\xi,\xi')\in \mathcal{O}_\varepsilon$. 
Equivalently
\begin{equation}
\label{sis1}
N\phi(x+\xi)=-\phi'(p+\xi').
\end{equation}

By definition of quasi-special, $\phi$ is Gauss-reduced of rank $r$.
Let 
\begin{equation*}
\phi=  \left(\begin{array}{c}
\phi_1\\
\vdots\\
\phi_r
\end{array}\right)=
 \left(\begin{array}{cccccc}
a&\dots &0&L_1\\
\vdots & &\vdots & \vdots \\
0&\dots &a&L_r
\end{array}\right)
\end{equation*}
with $L_i\in\mathbb{Z}^{g-r}$ and $H(\phi)=a$.

Since $\Gamma_p$ is the division group of $p$, the point  $y'$  defined as 
$$
N[a]y'=\phi'(p),$$ belongs to $\Gamma_p^r$.

Let $\upxi'$ be a point of $E^r$ such that 
$$
N[a]\upxi'=\tilde\phi(\xi,\xi').
$$
We define
\begin{equation*}\begin{split}
  y&=(y',0,\dots,0)\in \Gamma_p^r\times \{0\}^{g-r},\\
\upxi&=(\upxi',0,\dots,0)\in E^r\times \{0\}^{g-r}.
\end{split}\end{equation*}
We have 
\begin{equation*}\begin{split}
N\phi(y)&=N[a]y'=\phi'(p)\\
N\phi(\upxi)&=N[a]\upxi'=\tilde\phi(\xi,\xi').
\end{split}\end{equation*}
It follows
\begin{equation*}
\begin{split}
N\phi(x+y+\upxi)&=
N\phi(x)+\phi'(p)+\tilde\phi(\xi,\xi')\\ 
&=\tilde\phi((x,p)+(\xi,\xi'))=0.
\end{split}
\end{equation*}
Thus 
$$x\in C(\overline{\qe})\cap (B_{N\phi}+\Gamma^g_p+\mathcal{O}_{||\upxi||}).$$

In order to finish the proof, we shall prove $$||\upxi||\le \varepsilon K_4.$$
By definition of $ \upxi$ we see that

\begin{equation*}
\begin{split}
||\upxi||=||\upxi'||=\frac{||\tilde\phi(\xi,\xi')||}{Na} &\le (g+s)\frac{\max(H(\phi'),Na)}{Na}||(\xi,\xi')||\\
&\le (g+s)\frac{\max(H(\phi'),Na)}{Na} \varepsilon.
\end{split}
\end{equation*}
We claim \begin{equation*}
\frac{\max(H(\phi'),Na)}{Na}\le \frac{K_4}{g+s}.
\end{equation*}
Let $\phi'=(b_{ij})$. We shall prove that $H(\phi')=\max_{ij} |b_{ij}|\le \frac{K_4}{g+s} N a$.
Let $|b_{kl}|=H(\phi')$. 
Consider the $k$-row of the system (\ref{sis1}) 
$$N\phi_k(x)+N{\phi_k(\xi})=- {\sum_jb_{kj}(p_j+\xi'_j)}.$$
The triangle inequality gives
\begin{equation}
\label{icsk}
\frac{||\phi_k(x)||}{a}+\frac{||{\phi_k(\xi)}||}{a}\ge\frac {||\sum_jb_{kj}(p_j+\xi'_j)||}{Na}.
\end{equation}
Since $\varepsilon \le \eduep$ and $r\ge 2$, then $(x,p) \in  S_2(C \times p,\mathcal{O}_{\eduep})$ which has norm $\kduep$. Hence
$$||x||\le||(x,p)||\le \kduep.$$
Since $a=H(\phi)$, we see that 
$$\frac{||\phi_k(x)||}{a}\le(g-r+1)\kduep\,\,\,\,\,\,{\rm{and}}\,\,\,\,\,\, \frac{||{\phi_k(\xi)}||}{a}\le (g-r+1)\varepsilon.$$
Substituting in (\ref{icsk})
$$(g-r+1)(\kduep +\varepsilon)\ge \frac {||\sum_jb_{kj}(p_j+\xi'_j)||}{Na}.$$
Recall that $ \varepsilon \le \ezerop$. Hence, Proposition \ref{ret} with $(b_1,\dots,b_s)=(b_{k1},\dots,b_{ks})$ $(\xi_1, \dots ,\xi_s)=-\xi'$ and $\upxi=0$, implies
\begin{equation*}
\begin{split}(g-r+1)(\kduep +\varepsilon)&\ge\frac{1}{Na}\left( \cunop^2\sum_j|b_{kj}|^2||p_j||^2\right)^{\frac{1}{2}}
\\
&\ge\frac{{\cunop}H(\phi')}{Na}\min_i||p_i||.
\end{split}\end{equation*}
Whence $$H(\phi')\le \frac{K_4}{g+s} Na.$$

\end{proof}

The inclusion in Theorem \ref{equi34} ii. is proven only for  a
set $S_r(C\times p, \mathcal{O}_{\varepsilon})$ which is known to have
bounded height. 
If the norm $\kduep$   of $S_r(C\times p,\mathcal{O}_{\varepsilon})$
 goes  to infinity, the set 
$(\Gamma^g_p)_{\varepsilon\ktre}$ tends to be the whole $E^g$.

\begin{remark}
\label{optimal}
We would like to show that our Theorems  \ref{altzero} and \ref{tutto} are optimal. 
Let $\Gamma=\langle (y_1,0, \dots , 0)\rangle$,  where  $y_1$ is a non torsion point in $E(\overline\qe)$.
Since $C$ is transverse, the projection $\pi_1$ of $C(\overline\qe)$ on the first factor $E(\overline{\qe})$  is surjective. Let $x_n\in C(\overline\qe)$ such that  $\pi_1(x_n)=ny_1$. So $x_n-n(y_1,0,  \dots ,0)$ has first coordinate zero, and belongs to the algebraic subgroup $0\times E^{g-1}$.Then, for all $n\in \mathbb{N}$ it holds
$$x_n\in B_{\phi=(1,0,\dots,0)}+\Gamma.$$
This shows that $x_n \in S_1(C,\Gamma)$, so $S_1(C,\Gamma)$ does not have bounded height.
By Theorem \ref{equi34} part i, neither $S_1(C\times y_1)$ has bounded height.

\end{remark}

\section{  Special Morphisms and an important Inclusion}
We can actually show a  stronger inclusion than the one in Theorem
\ref{equi34} i. The set $S_r(C,\Gamma_\varepsilon)$ can be included in
a  subset of $S_r(C\times \gamma, \mathcal{O}_\varepsilon)$, namely the
subset defined by special morphisms.

\begin{D}[Special Morphisms]
\label{spec}
A surjective morphism  $\tilde{\phi}:E^{g+s} \to E^r$ is special if $\tilde\phi=(N\phi|\phi')$ is quasi-special and satisfies the further condition
$$H(\tilde\phi)=NH(\phi).$$
Equivalently $\tilde\phi$ is special if and only if
\begin{enumerate}
\item
$\tilde\phi$ is Gauss-reduced,
\item $ H(\tilde\phi)I_r$  is a submatrix of the matrix consisting of the first $g$ columns of $\tilde\phi$.
\end{enumerate}

\end{D}
Let us prove the equivalence of these two definitions:
\begin{proof}
That the first definition implies the second  is a clear matter.
For the converse,
decompose $\tilde\phi=(A|\phi')$ with $A\in M_{r\times g}(\ze)$ and
$\phi'\in M_{r\times s}(\ze)$. Let $N$ be the greatest common divisor
of the entries of $A$. Define $\phi=A/N$ and $a=H(\tilde\phi)/N$. Then $\phi=(aI_r|L')$ 
 is Gauss-reduced and $\tilde\phi=(N\phi|\phi')$.
 \end{proof}

A nice remark is that the obstruction to show unconditionally that $S_r(C\times p,\mathcal{O}_\varepsilon)$ is included in $S_r(C,(\Gamma_p^g)_{\varepsilon'})$ is exactly due to the non-special morphisms.  Sets of the kind
$$\left(C(\overline{\qe})\times p\right)\cap
\left(B_{\tilde\phi}+\mathcal{O}_\varepsilon\right)$$ which do not have bounded
height, can be included in
$S_r(C,(\Gamma_p^g)_{\varepsilon'})$  if $\tilde\phi$ is special, indeed in general $\varepsilon'={\rm{c}}(g,s)\frac{H(\tilde\phi)}{H(A)}\varepsilon$ for any $\tilde\phi=(A|\phi')$.

\begin{propo}
\label{speciali}
Let $2\le r$ and  $\varepsilon\le  \min (\euno,\frac{\kuno}{g})$. The map $x \to (x,\gamma)$ defines  an injection 
$$\bigcup_{{\substack{\phi\,\,\,{\rm{Gauss \,\,\,reduced}}\\ {\rm{rk}}( \phi)=r}} } C(\overline{\qe})\cap \left(B_\phi+\left(\Gamma^g_0\right)_{\varepsilon}\right)
\hookrightarrow 
\bigcup_{\substack{\tilde\phi=(N\phi|\phi')\,\,\, {\rm{special}}\\ {\rm{rk}}\tilde\phi=r }}
(C(\overline{\qe})\times \gamma) \cap (B_{\tilde{\phi}} +\mathcal{O}_{\varepsilon}).$$

\end{propo}
\begin{proof}

  Let $x \in C(\overline{\qe})\cap (B_\phi+\Gamma^g_0 +\mathcal{O}_{\varepsilon})$ with $\phi$ Gauss-reduced of rank $r$. Then, there exist $y\in \Gamma^g_0$ and $\xi \in \mathcal{O}_{\varepsilon}\subset E^g$ such that 
$$\phi(x+y+\xi)=0.$$
Since $\gamma_1, \dots ,\gamma_s$ is a maximal free set of $\Gamma_0$, there exists an integer $N$ and a matrix $G \in M_{r,s}(\mathbb{Z})$ such that 
$$[N]y=G(\gamma).$$
Let $n$ be the greatest common divisor of the entries of $({N}\phi | \phi  G)$.
We define$$\tilde{\phi}=\frac{1}{n}({N}\phi | \phi  G).$$
Clearly 
\begin{equation*}
\begin{split}
({N}\phi | \phi  G)\left((x,\gamma)+(\xi,0)\right)&=
N\phi(x)+\phi G(\gamma)+N\phi(\xi)\\
&=N\phi(x+y+\xi)=0.
\end{split}
\end{equation*}
Thus \begin{equation}
\label{real1}
 n\tilde\phi((x,\gamma)+(\xi,0))=0.
 \end{equation}
Equivalently
$$(x,\gamma)\in  (C(\overline{\qe})\times \gamma) \cap (B_{n\tilde{\phi}} +\mathcal{O}_{\varepsilon}).$$
By Lemma \ref{tor} i.  with $\psi=\tilde\phi$ and $N=n$, it follows
$$(x,\gamma)\in  (C(\overline{\qe})\times \gamma) \cap (B_{\tilde{\phi}} +\mathcal{O}_{\varepsilon}).$$

We shall still  show  that  $\tilde{\phi}$ is special, using the first definition of special. 
By assumption, the morphism $\phi$ is Gauss-reduced.
By definition of $\tilde\phi$, the greatest common divisor of its
entries is $1$.
In order to conclude that $\tilde\phi$ is special, we  still have to show that $$H(\tilde\phi)=N a$$ or equivalently $$H(\phi')\le N a.$$

The proof is similar to the last part of the proof of Theorem \ref{equi34} ii.

Let $\phi'=(b_{ij})=\phi G$. Let $|b_{kl}|=\max_{ij}
|b_{ij}|=H(\phi')$. Let $\phi_k$ be the $k$-th
row of $\phi$. Consider the $k$-th row of the system (\ref{real1})
\begin{equation*}
nN(\phi_k(x)+{\phi_k(\xi)} )=-n {\sum_jb_{kj}\gamma_j}.
\end{equation*}
Then $$\frac{||\phi_k(x)||}{a}+\frac{||\phi_k(\xi)||}{a}\ge \frac{1}{Na}{\Big|\Big|\sum_jb_{kj}\gamma_j\Big|\Big|}.$$
Clearly $x \in S_r(C, (\Gamma_0^g)_\varepsilon)$. Since $\varepsilon \le \euno$, then  $x \in S_2(C,(\Gamma_0^g)_{\euno})$ which has norm bounded by $\kuno$. So $$|| x||\le \kuno.$$ 
As $H(\phi_k)\le H(\phi)=a$, $$\frac{||\phi_k(x)||}{a}\le (g-r+1)K_1.$$ Furthermore
$$\frac{||\phi_k(\xi)||}{a}\le (g-r+1)\varepsilon.$$
Then
 $$(g-1)(\kuno+\varepsilon)\ge \frac{1}{Na}{\Big|\Big|\sum_jb_{kj}\gamma_j\Big|\Big|}.$$

From relations (\ref{con1}) with $(b_1,\dots ,b_s)=(b_{k1},\dots ,b_{ks})$ and (\ref{bignor}) in the recap, we deduce
\begin{equation*}
\begin{split}
(g-1)(\kuno+\varepsilon)&\ge \frac{1}{Na}\left(\frac{1}{9} \sum_j|b_{kj}|^2||\gamma_j||^2\right)^{\frac{1}{2}}
\\&\ge\frac{ H(\phi')}{3Na}\min_j||\gamma_j||\\&\ge \frac{ H(\phi')}{3Na}{3gK_1}.
\end{split}\end{equation*}
We  assumed $\varepsilon\le \frac{K_1}{g}$, so $$H(\phi')\le Na.$$

\end{proof}

This inclusion is important; the Bogomolov type bounds are given for intersections with
$\mathcal{O}_\varepsilon$ and not with $\Gamma_\varepsilon$. Actually
there exist bounds for $\varepsilon$, such that  $C\cap \Gamma_\varepsilon$ is finite. 
These bounds are deduced by the Bogomolov type bounds and their
dependence on the degree of the curve is  not sharp
enough for our purpose. 
 To overcome such an obstacle and
solve the problem with $\Gamma_\varepsilon$, we use the above Proposition
\ref{speciali} and make use of the Bogomolov type bounds for $C\times
\gamma$ intersected with $B_{\tilde{\phi}} +\mathcal{O}_{\varepsilon}$, where $\tilde \phi $ is special of rank $2$.

\section{The Proof of  Theorem \ref{tutto}}
\label{cg}

In  sections \ref{subunion} and \ref{due} below, we prepare the core for the proof of
 Theorem \ref{tutto}. 
In Proposition $\centro$   we prove that the  union can be taken over    finitely
many sets, while in Proposition  $\finito$ we prove that each set in the union is finite. Hence, our set is  finite.

We prefer to present first the proof of Theorem \ref{tutto}  and then to prove the two key Propositions $\centro$ and $\finito$. We hope that,  knowing a priory the  aim of sections  \ref{subunion} and \ref{due},  the reader gets the right inspiration to handle them.

\begin{proof}[Proof of Theorem \ref{tutto}]
Assuming Conjecture \ref{sin2}, we prove Conjecture \ref{geni} part iv. In view of  Theorem \ref{equivalenza}, also   part iii. is proven. Parts i. and ii. are then  obtained by setting $\varepsilon=0$.

Choose

\begin{enumerate}
 \item  $n=2(g+s)-3$.
 \item 
$ \delta_1 =\frac{\min(\ecinque,\edue)}{(g+s)^2} $ \hfill
  where  
  $ \ecinque\,\,\,{\rm{is \,\,\, as\,\,\,in\,\,\,Proposition \,\,\,}}\finito,$ \\
  \item $\delta= {\delta_1}{{M'}^{-1-\frac{1}{2n}}}$  \hfill{where
 $M'=\max\left(2,\lceil \frac{\kdue}{\delta_1}\rceil^2\right)^{n}.$}

 \end{enumerate}
\vspace{0.5 cm}
Since $\Gamma_\delta\subset (\Gamma_0^g)_\delta$, 
$$S_2(C,\Gamma_{\delta})\subset S_2(C,(\Gamma_0^g)_{\delta}).$$
In Lemma \ref{tor1} with $\varepsilon=\delta$, we saw that 
$$S_2(C,(\Gamma_0^g)_{\delta}) \subset \bigcup_{\substack{\phi\,\,\,{\rm{Gauss-reduced}}\\ {{\rm{rk}}}\phi=2}}  C(\overline{\qe})\cap (B_\phi+(\Gamma_0^g)_{\delta}).$$
Note that   $\delta<\delta_1< \min(\euno,\frac{\kuno}{g})$. Then, Proposition \ref{speciali} with $\varepsilon=\delta$ implies  that 
$$\bigcup_{\substack{\phi\,\,\,{\rm{Gauss-reduced}}\\ {{\rm{rk}}}\phi=2}}  C(\overline{\qe})\cap (B_\phi+(\Gamma_0^g)_{\delta}) \hookrightarrow  \bigcup_{\substack{\tilde\phi=(N\phi|\phi')\,\,\, {\rm{special}}\\\,\,\,{{\rm{rk}} }\phi=2 }}(C(\overline{\qe})
 \times \gamma) \cap (B_{\tilde{\phi}}
  +\mathcal{O}_{\delta}).$$
Note that $\delta= {\delta_1}{{M'}^{-(\expalt)}}$ and $\delta_1\le \edue$. Then,  Proposition $\centro$    ii.  in section \ref{subunion} below,  with $\varepsilon=\delta_1$, $r=2$ (and $n$ is already defined as $2(g+s)-4+1$), shows that 
$$ \bigcup_{\substack{\tilde\phi\,\,\, {\rm{special}}\\ \,\,\,{{\rm{rk}} }\phi=2 }}(C(\overline{\qe})
 \times \gamma) \cap (B_{\tilde{\phi}}
  +\mathcal{O}_{\delta})$$ is a subset of 
  \begin{equation}
  \label{uno}
  \bigcup_{\substack{\tilde\phi \,\,\, {\rm{special}}\\H(\tilde\phi)\le M'\,\,\,{{\rm{rk}} }\phi=2 }}(C(\overline{\qe})
 \times \gamma) \cap \left(B_{\tilde{\phi}}
  +\mathcal{O}_{(g+s)\delta_1/H(\tilde\phi)^\expalt}\right).
  \end{equation}
 Observe that  in (\ref{uno}), $\tilde\phi$ ranges over finitely many morphisms, as $H(\tilde\phi)$ is bounded by $M'$.
  
 We have chosen $\delta_1 \le \frac{\ecinque}{(g+s)^2}$. Proposition $\finito$   ii.   in  section \ref{due} below with $\varepsilon=(g+s)\delta_1$, implies that  for all $\tilde\phi=(N\phi|\phi')$ special of rank $2$, the set $$(C(\overline{\qe})
 \times \gamma) \cap \left(B_{\tilde{\phi}}
  +\mathcal{O}_{(g+s)\delta_1/H(\phi)^\expalt}\right)$$ is finite.  Note that $H(\phi)\le H(\tilde\phi)$, thus also the sets $$(C(\overline{\qe})
 \times \gamma) \cap \left(B_{\tilde{\phi}}
  +\mathcal{O}_{(g+s)\delta_1/H(\tilde\phi)^\expalt}\right)$$ appearing in (\ref{uno}) are  finite.
 
It follows that,
 the set $S_2(C,\Gamma_\delta)$ is contained in the union of finitely many
finite sets. 
 So it is finite.

\end{proof}

Despite our proof relying on  Dirichlet's
Theorem and a Bogomolov type bound, a direct use of these two theorems is not sufficient to prove
Theorem \ref{tutto}. 
Using   Dirichlet's Theorem in a more natural way, one can prove that, for $r\ge2$,
$$S_r(C,\Gamma_\varepsilon)\subset \bigcup_{H(\phi)\le M(\varepsilon),\,\,\,{\rm{rk}}\phi=r}C(\overline{\qe})\cap \left(B_\phi+\Gamma_\varepsilon\right).$$
On the other hand, a  direct use of Bogomolov's type bound gives that  
$$C(\overline{\qe})\cap \left(B_\phi+\mathcal{O}_{\varepsilon/H(\phi)^2}\right)$$ is finite, for $\phi$ of rank at least $2$.
Even if we forget $\Gamma$, the discrepancy between $\varepsilon$ and
$\varepsilon/H(\phi)^2$ does  not look encouraging, and it took us a
long struggle to overcome the problem. In  Propositions $\centro$ and $\finito$, we succeed in overcoming the mismatch; in both statements we  obtain neighbourhoods of radius $\varepsilon/H(\phi)^\expalt$.

Do not be misled by the  following wrong thought: One might
think that, since we consider only morphisms $\phi$ such that
$H(\phi)\le M$, it could be enough to choose $\varepsilon'=\varepsilon/M^2$. However, $M=M(\varepsilon)$ is an unbounded function of $\varepsilon$ as  $\varepsilon$ tends to $0$.

\section{PART I: The Box Principle and the Reduction to a finite Sub-Union} 
\label{subunion}
 In Lemma \ref{dicov}, we   approximate a Gauss-reduced morphism with a Gauss-reduced
morphism of bounded height. On a set of
bounded height,
such an approximation allows us to consider unions over finitely many algebraic subgroups, instead of unions over all algebraic subgroups (see Proposition $\centro$ below).

We recall Dirichlet's Theorem on the rational approximation of reals. 

\begin{thm1}[Dirichlet 1842, see  \cite{S} Theorem 1 p. 24]

\label{diry}
Suppose that $\alpha_1, \dots , \alpha_n$ are $n$ real numbers and that $Q\ge 2$ is an integer. Then there exist integers $f, f_1,\dots,f_n$ with
\begin{equation*}
1\le f<Q^n \,\,\,\,{\rm{and}}\,\,\,\,\left|\alpha_if-f_i\right|\le \frac{1}{Q}
\end{equation*}
for $1\le i \le n$.

\end{thm1}

\begin{lem}
\label{dicov}
Let $Q\ge 2$ be an integer. Let $\fiffo= (aI_r|L)\in M_{r\times g}(\ze)$ be  Gauss-reduced. 
Then, there exists a Gauss-reduced   $\pippo=(fI_r|L')\in M_{r\times g}(\ze)$ such that
\begin{enumerate}
\item $H(\pippo)=f \le Q^{rg-r^2+1},$
\item $\left| \frac{\pippo}{f}- \frac{\fiffo}{a} \right|\le {Q^{-\frac{1}{2}}f^{-1-\frac{1}{2(rg-r^2+1)}}}.$
\end{enumerate}
The norm  $|\cdot|$  of a matrix  is the maximum of the
absolute values of its entries.

\end{lem}
\begin{proof}
If $a\le Q^{rg-r^2+1}$, no approximation is needed  as $\fiffo$ itself
satisfies the conclusion. So we can assume that
\begin{equation*}
\fiffo=
 \left(\begin{array}{cccccc}
a&\dots &0&L_1\\
\vdots & & \vdots & \vdots\\
0&\dots &a&L_r
\end{array}\right)
\end{equation*}
  is a Gauss-reduced morphism such that  $H(\fiffo)=a> Q^{rg-r^2+1}$.
  Consider the element  $$\alpha= \left(1,\frac{L_1}{a},\dots ,
    \frac{L_r}{a}\right)=(\alpha_1,\alpha_2, \dots , \alpha_{rg-r^2+1})
 $$
  in  $\mathbb{R}^{rg-r^2+1}$. 

 Define $n=rg-r^2+1$. Apply Dirichlet's Theorem to $\alpha$. Then, there
  exist integers $f, f_1,\dots,f_n$ with
\begin{equation}
\label{dir}
1\le f<Q^n \,\,\,\,{\rm{and}}\,\,\,\,\left|\alpha_if-f_i\right|\le \frac{1}{Q}
\end{equation}
for $1\le i \le n$.
We can assume that $f, f_1, \dots, f_n$ have greatest common divisor $1$.
Define $$w=\frac{1}{f}(f_1, \dots , f_n)=\frac{1}{f}(f_1,L'_1, \dots,L'_r),$$ with $L'_i\in \mathbb{Z}^{g-r}$.
We claim that
\begin{equation*}
\begin{split}
f_1&=  f,\\ 
|f_i|&\le f.
  \end{split}
  \end{equation*} In fact, by (\ref{dir}) for $i=1$ we have
$\Big|\frac{f_1}{f}- 1\Big|\le \frac{1}{Qf}$, which implies that $|f- f_1|<1$. Since $f$ and $f_1$ are integers, we must have 
$f=  f_1$. Similarly, by (\ref{dir}) for $i=2,\dots,n$, we have $\left|\frac{f_i}{f}- \alpha_i\right|
\le\frac{1}{Qf}$. This implies that $|f_i|\le f+\frac{1}{Q}$.  We deduce $|f_i|\le f$.

It follows that 
\begin{equation*}
\pippo= 
 \left(\begin{array}{cccc}
f&\dots &0&L'_1\\
\vdots & & \vdots \\
0&\dots &f&L'_r
\end{array}\right)
\end{equation*}
is a Gauss-reduced morphism of rank $r$ with $H(\pippo)=f$. 

Relation (\ref{dir}) immediately gives 
$$f\le Q^{n}$$
and 
$$\left|\frac{\pippo}{f}-\frac{\fiffo}{a}\right| \le \frac{1}{Qf}\le \frac{1}{Q^{\frac{1}{2}}f^{1+\frac{1}{2n}}},$$
where in the last inequality we use $Q^{\frac{1}{2}}\ge f^{\frac{1}{2n}}$.

\end{proof}

At last we prove our first main proposition; the union can be taken over finitely many algebraic subgroups.

 If $\fiffo$ has large height and $B_\fiffo$ is close to  $x$, with $x$ in  a set of bounded height, then there
 exists $\pippo $  with height bounded by a constant  such that $B_\pippo$   is also
 close to $x$. One shall be careful that, in the following inclusions,
 on the left hand side we consider a neighbourhood of $B_\fiffo$ of
 fixed radius, while on the right hand side the neighbourhood becomes
 smaller as the height of $\pippo$ grows. This is a crucial gain, with
 respect to the simpler approximation (obtained by a direct use of
 Dirichlet's Theorem) where the neighbourhoods have constant radius on both hand sides.

\begin{propb*} 
\label{centro}

 Assume $r\ge 2$.
 \begin{enumerate}

\item
If $0<\varepsilon\le  \euno$, then

$$\bigcup_{ \fiffo}   C(\overline{\qe})\cap \left(B_\fiffo+\left(\Gamma^g_0\right)_{\varepsilon/M^\expaltu}\right) \subset \bigcup_{{H(\pippo)\le M} }  C(\overline{\qe})\cap \left(B_\pippo+\left(\Gamma^g_0\right)_{g\varepsilon/H(\pippo)^\expaltu}\right) ,$$

\noindent where $\phi$ and $\psi$ range over Gauss-reduced morphisms of rank $r$, $n=rg-r^2+1$ and $M= \max\left(2,\lceil\frac{\kuno}{\varepsilon}\rceil^2\right)^{n}$. 

\vspace{0.4cm}

\item
 If $0<\varepsilon\le  \edue$, then

$$ \bigcup_{\tilde\fiffo}(C(\overline{\qe})\times \gamma) \cap \left(B_{\tilde{\fiffo}} +\mathcal{O}_{\varepsilon/ {M'}^\expaltu}\right) \subset 
 \bigcup_{H(\tilde\pippo)\le M'}(C(\overline{\qe})
 \times \gamma) \cap \left(B_{{\tilde\pippo}}
  +\mathcal{O}_{(g+s)\varepsilon/H(\tilde\pippo)^\expaltu}\right) ,$$
  where  $\tilde\phi$ and $\tilde\psi$ range  over special morphisms of rank $r$, $n=r(g+s)-r^2+1$ and $M'=\max\left(2,\lceil\frac{\kdue}{\varepsilon}\rceil^2\right)^{n}$.

\end{enumerate}

\end{propb*}

\begin{proof}

{\bf Part i.}  Let $\fiffo=(aI_r|L)$ be Gauss-reduced of rank $r$.

First consider the case $H(\fiffo) \le M$. Then ${\varepsilon/M^\expaltu}\le {\varepsilon/H(\fiffo)^\expaltu}$. Obviously $$C(\overline{\qe})\cap (B_\fiffo+\Gamma^g_0 +\mathcal{O}_{\varepsilon/M^\expaltu})\subset C(\overline{\qe})\cap (B_\fiffo+\Gamma^g_0 +\mathcal{O}_{\varepsilon/H(\fiffo)^\expaltu})$$ is contained in the right hand side.

Secondly consider the case $H(\fiffo) > M$. We shall show that  there exists $\pippo$ Gauss-reduced with $H(\pippo)\le M$ such that 
$$C(\overline{\qe})\cap \left(B_\fiffo+\Gamma^g_0 +\mathcal{O}_{\varepsilon/M^\expaltu}\right)\subset C(\overline{\qe})\cap \left(B_\pippo+\Gamma^g_0 +\mathcal{O}_{g\varepsilon/H(\pippo)^\expaltu}\right).$$

We fix $Q=\max\left(2,\lceil\frac{\kuno}{\varepsilon}\rceil^2\right)$. Recall that $n=rg-r^2+1$.
By Lemma \ref{dicov}, there exists 
a Gauss-reduced morphism 
\begin{equation*}
\pippo=
 \left(\begin{array}{cccccc}
f&\dots &0&L'_1\\
\vdots & & \vdots \\
0&\dots &f&L'_r
\end{array}\right)
\end{equation*}
such that 
$$H(\pippo)=f\le M$$ 
and 
\begin{equation}
\label{borne}
\left|\frac{\pippo}{f}-\frac{\fiffo}{a}\right|\le \frac{1}{Q^{\frac{1}{2}}f^{1+\frac{1}{2n}}}.
\end{equation}

Let $x \in  C(\overline{\qe})\cap(B_\fiffo+\Gamma^g_0+ \mathcal{O}_{{\varepsilon}/M^\expaltu})$. Then there  exist $y\in \Gamma^g_0$ and $\xi\in \mathcal{O}_{{\varepsilon}/M^\expaltu}$ such that 
\begin{equation*}
\fiffo(x-y-\xi)=0.
\end{equation*}
We want to show that there exist  $y' \in {\Gamma^g_0}$ and
$\xi' \in  \mathcal{O}_{{g\varepsilon/f^\expaltu}}$ such that 
$$\pippo(x-y'-\xi')=0.$$

Let $y'' $ be a point such that  $$[a] y'' ={\fiffo (y)}.$$ As $\Gamma_0$ is a division group,
$y'' \in {\Gamma^r_0}$. We define $$y'=(y'',0)\in  {\Gamma^r_0}\times\{0\}^{g-r},$$ then
 \begin{equation*}
\pippo(y')=[f]y''.
\end{equation*}

Let $\xi'' $ be a point such that 
$$[f]\xi'' ={\pippo (x-y')}.$$ 
We define $$\xi'=(\xi'',0),$$ then
$$\pippo(\xi')=[f]\xi''=\pippo(x-y').$$
So
$$
\pippo(x-y'-\xi')=0.$$
It follows   that
$$ x \in C(\overline{\qe})\cap (B_\pippo+\Gamma^g_0 +\mathcal{O}_{||\xi'||}).$$
In order to finish the proof, we are going to prove that 
$$||\xi'||\le\frac{g\varepsilon}{f^\expaltu}.$$
By definition
\begin{equation*}
||\xi'||=||\xi''||=\frac{||\pippo (x-y')||}{f}.
\end{equation*}
Consider the equivalence
 \begin{equation*}
\begin{split}
a\pippo (x-y') & = a\pippo(x)-a\pippo(y')\\ 
& = a\pippo(x)-a[f]y''\\
&=a\pippo(x)-f\fiffo(y)\\
&=a{\pippo (x)}-f{\fiffo (x)}+ f\fiffo(\xi).
\end{split}
\end{equation*}
Then
\begin{equation*}
||\xi'||=\frac{1}{af}||f\fiffo(\xi)-f{\fiffo (x)}+a{\pippo (x)}||\le \frac{1}{a}||\fiffo(\xi)||+\frac{1}{af}||a\pippo(x)-f{\fiffo (x)}||.
\end{equation*}
Let us estimate separately each norm on the right.

On one hand $$\frac{1}{a}||\fiffo(\xi)||\le (g-r+1)||\xi||\le \frac{(g-1)\varepsilon}{M^\expaltu}\le \frac{(g-1)\varepsilon}{f^\expaltu},$$
because $||\xi||\le \varepsilon/M^\expaltu$ and $f\le  M.$

On the other hand, since the rank of $\fiffo$ is at least 2 and $\varepsilon\le \euno$, we have that  $x\in S_2(C,(\Gamma^g_0)_{\euno})$, which has norm $K_1$. Thus
\begin{equation*}
||x||\le \kuno.
\end{equation*}
Using relation (\ref{borne}) and that $Q\ge\lceil \frac{\kuno}{\varepsilon}\rceil^2$, it follows that
\begin{equation*}
\begin{split}
\frac{1}{af}||a{\pippo (x)}-f{\fiffo (x)}||&\le
\left|\frac{\pippo }{f}-\frac{\fiffo }{a}\right|||x||\\ 
&\le \frac{1}{Q^{\frac{1}{2}}f^{\expaltu}}||x||\\
&\le \frac{\varepsilon||x||}{\kuno f^{\expaltu}}\\
&\le  \frac{\varepsilon}{f^{\expaltu}}.
\end{split}
\end{equation*}

We conclude that 
$$||\xi'||\le \frac{(g-1)\varepsilon}{f^\expaltu}+\frac{\varepsilon}{f^{\expaltu}}
\le\frac{g\varepsilon}{f^\expaltu}.$$

 {\bf Part ii.} 
We fix $Q=\max(2,\lceil\frac{\kdue}{\varepsilon}\rceil^2)$. 

Let $\tilde{\fiffo}=(N\fiffo|\fiffo'):E^{g+s} \to E^r$ be special. By  the second definition of special
$$\tilde{\fiffo}=(\enne I_r|\ast)$$ is Gauss-reduced and $H(\tilde\fiffo)=\enne$.

As in part i.; if $H(\tilde\fiffo)\le M'$ then $\varepsilon/{M'}^\expaltu\le \varepsilon/H(\tilde\fiffo)^\expaltu$ and 
$$ (C(\overline{\qe})\times \gamma)\cap \left(B_{\tilde{\fiffo}}+\mathcal{O}_{\varepsilon/{M'}^\expaltu}\right)$$ is contained in the right hand side.

 Now, suppose  that $H(\tilde\fiffo)> M'$. Recall that $n=r(g+s)-r^2+1$.
 By
 Lemma \ref{dicov} (applied with $\fiffo=\tilde\fiffo$ and $\pippo=\tilde\pippo$)
there exists
$\tilde{\pippo}=(fI_r|\ast)$  Gauss reduced such that  
$$H(\tilde\pippo)=f\le M'$$ 
and 
\begin{equation}
\label{rel1}
\left|\frac{\tilde{\fiffo}}{\enne}-\frac{\tilde{\pippo}}{f}\right|\le \frac{1}{Q^{\frac{1}{2}}f^{1+\frac{1}{2n}}}.
\end{equation}
Then $\tilde\pippo$ is special, according to the second formulation in Definition \ref{spec}.

The proof is now similar to the proof of  part i. 
 We want to show that, if $${\tilde{\fiffo}}((x,\gamma) +\xi)=0$$ for $\xi \in\mathcal{O}_{\varepsilon/{M'}^\expaltu}$, then $${\tilde{\pippo}}((x,\gamma) +\xi')=0$$ for $\xi'\in \mathcal{O}_
  {(g+s)\varepsilon/H(\tilde\pippo)^\expaltu}$.
 
Let  $\xi''$ be a point in $ E^r$ such that $$
[f]\xi''=-\tilde\pippo(x,\gamma).$$
Let $\xi'=(\xi'',\{0\}^{g-r+s})$. Then
$$\tilde\pippo((x,\gamma)+(\xi',0))=0.$$ It follows
$$(x,\gamma)\in (C(\overline{\qe})\times \gamma)\cap (B_{\tilde{\pippo}}+\mathcal{O}_{||\xi'||}),$$
where $\tilde\pippo$ is special and $H(\tilde\pippo)\le M'$.

It remains to prove that
  $$||\xi'|| \le\frac{(g+s)\varepsilon}{H(\tilde\pippo)^\expaltu}.$$  
Obviously
\begin{equation*}
\enne\tilde\pippo(x,\gamma)  = \ti\left({\tilde\fiffo(x,\gamma)}-{\tilde\fiffo(x,\gamma)}\right)+\enne{\tilde\pippo(x,\gamma)}.
\end{equation*}
According to the definition of $\xi'$,
\begin{equation*}
\begin{split}||\xi'||=||\xi''||=\frac{||\tilde\pippo(x,\gamma)||}{f}&=\frac{1}{\enne\ti}\left|\left|\ti\left({\tilde\fiffo(x,\gamma)}-{\tilde\fiffo(x,\gamma)}\right)+\enne{\tilde\pippo(x,\gamma)}\right|\right|\\
&\le \frac{1}{\enne}\Big|\Big|{\tilde\fiffo(x,\gamma)}\Big|\Big|+\frac{1}{\enne\ti}\Big|\Big|\enne{\tilde\pippo(x,\gamma)}-f{\tilde\fiffo(x,\gamma)}\Big|\Big|.
\end{split}
\end{equation*}
We estimate the two norms on the right. 

On one hand
\begin{equation*}
\begin{split}\frac{||\tilde\fiffo(x,\gamma)||}{\enne}= \frac{||\tilde\fiffo(\xi)||}{\enne} &\le (g-r+1+s) ||\xi||\\
&\le \frac{(g-r+1+s)\varepsilon}{{M'}^\expaltu} \\&\le \frac{(g-r+1+s)\varepsilon}{\ti^\expaltu},
\end{split}
\end{equation*}
where in the last inequality we use that  $\ti\le M'$. 

On the other hand,
by definition of $\edue$,  we know that the norm of the set $S_2(C\times \gamma, \mathcal{O}_{\edue})$  is bounded by $\kdue$. 
 Since $\varepsilon \le \edue$,  we have $(x,\gamma)\in S_2(C\times \gamma,
 \mathcal{O}_{\edue})$. Therefore $$||(x,\gamma)||\le \kdue.$$

Using relation (\ref{rel1}) and that $Q\ge \lceil\frac{\kdue}{\varepsilon}\rceil^2$, we estimate
\begin{equation*}
\begin{split}
\frac{1}{\enne\ti}\Big|\Big|\enne{\tilde\pippo(x,\gamma)}-\ti\tilde\fiffo(x,\gamma)\Big|\Big|&\le
\left|\frac{\tilde\fiffo}{\enne}-\frac{\tilde\pippo}{\ti}\right|||(x,\gamma)||\\&\le
\frac{||(x,\gamma)||}{{Q^{\frac{1}{2}}f^{1+\frac{1}{2n}}}}\\
&\le
\frac{\varepsilon ||(x,\gamma)||}{(\kdue) \ti^\expaltu}\le
\frac{\varepsilon}{\ti^\expaltu}.
\end{split}
\end{equation*}

Since $r\ge2$, we conclude

\begin{equation*}
||\xi'||\le \frac{(g-1+s)\varepsilon}{\ti^\expaltu}+\frac{\varepsilon}{\ti^\expaltu}=\frac{(g+s)\varepsilon}{H(\tilde\pippo)^\expaltu}.
\end{equation*}

\end{proof}

\section{PART II: The essential Minimum and  the Finiteness of each Intersection}
\label{due}

  Up until now we have used, several times, the fact  that the height of  our sets is
  bounded  (Theorem \ref{altzero}). In this section
  we often use that we work with a curve.
  
  In the following, we set  $$n=2(g+s)-3.$$
We would like to use Conjecture \ref{sin2} in order to provide $\varepsilon>0$ such that,  for all $\phi$ Gauss-reduced of rank $r=2$, the set 
 \begin{equation}
 \label{insi}
 \cqp\cap
\left(B_\phi+\mathcal{O}_{\varepsilon/H(\phi)^\expalt}\right)
\end{equation}
is finite.
This set  is simply 
$$\phi_{|\cp}^{-1}\left(\phi(\cp)\cap \phi\left(\mathcal{O}_{\varepsilon/H(\phi)^\expalt}\right)\right).$$
Further $$\phi\left(\mathcal{O}_{\varepsilon/H(\phi)^\expalt}\right)\subset
\mathcal{O}_{g\varepsilon/H(\phi)^\gain},$$  because if $\upxi \in
\mathcal{O}_{\varepsilon/H(\phi)^\expalt}$ then $||\phi(\upxi)||\le g H(\phi) ||\upxi||\le  g \varepsilon H(\phi)^{-\gain}$.
Thus, the set (\ref{insi})
is contained in the preimage 
of $$\phi(\cp)\cap \mathcal{O}_{ g \varepsilon/H(\phi)^\gain}.$$
If we can ensure that there exists $\varepsilon>0$ such that, for all morphisms $\phi$ Gauss-reduced of rank $r=2$,
\begin{equation}\label{ge} g  \varepsilon H(\phi)^{-\gain} < \mu(\phi(\cp)),\end{equation} then  the set 
 (\ref{insi}) is finite.

It is noteworthy   that  a direct use of a Bogomolov type bound, even optimal, is not
successful  in the following sense:  For a curve $X\subset E^g$ and any $\eta>0$, Conjecture \ref{sin2}  provides an invariant $\epsilon(X,\eta)$ such that $\epsilon(X,\eta)<\mu(X)$.
To ensure (\ref{ge}), we could naively require that 
\begin{equation*}
\label{geprimo}
g\varepsilon H(\phi)^{-\gain} \le \epsilon(\phi(\cp),\eta)\end{equation*}
 for all $\phi$ Gauss-reduced of rank $r=2$. Nevertheless 
this can be fulfilled only for $\varepsilon=0$.  


 

We need to throw new light on the problem to prove (\ref{ge}); via some isogenies, we construct a helping-curve $D$ and then we relate its
essential minimum to $\cp$. We then apply Conjecture \ref{sin2} to $D$. Thus,
we manage to provide a good lower bound for the essential minimum of $\cp$. We precisely take advantage of the fact that   $\mu([b]C)=b\mu(C)$, while $\epsilon([b]C,\eta)=\frac{\epsilon(C,\eta)}{b^{\frac{1}{g-1}+2\eta}}$ for any positive integers $b$.

\vspace{0.4cm}

Let 
\begin{equation*}
\phi= \left(\begin{array}{c}
\phi_1\\
\phi_2
\end{array}\right)= 
 \left(\begin{array}{ccccc}
a & 0& L_1\\
0 & a& L_2
\end{array}\right)
\end{equation*}
be a Gauss-reduced morphism of rank $2$ with $H(\phi)=a$.
We denote by   $\overline{x}=(x_3,\dots ,x_g)$, and recall that $n=2(g+s)-3$. 

We define
$$a_0=\lfloor a^\gain\rfloor.$$
 We associated to the morphism $\phi$ an isogeny
  \begin{equation*}
  \begin{split}
   \Phi:&E^g\to E^g\\
  &(x_1,\dots , x_g)\to (a_0\phi(x),x_3, \dots ,x_g).
    \end{split}
  \end{equation*}
We then relate it to  the isogenies:
  \begin{equation*}
  \begin{split}
  A: &E^{g} \to E^g\\
  & (x_1, \dots , x_g)\to(x_1,x_2,ax_3, \dots , ax_g).\\[0.3 cm]
   A_0: &E^{g} \to E^g\\
  & (x_1, \dots , x_g)\to(x_1,x_2,a_0x_3, \dots , a_0x_g).\\[0.3 cm]
   L:&E^g\to E^g\\
  &(x_1,\dots , x_g)\to (x_1+L_1(\overline{x}),x_2+L_2(\overline{x}),x_3, \dots ,x_g).    \end{split}
  \end{equation*}

\begin{D}[Helping-curve]
We define the curve $D$ to be  an  irreducible component of $A_0^{-1}LA^{-1}(C)$, where
$(\cdot)^{-1}$ simply means the pre-image.
\end{D}
The obvious relation
$$ [a_0a]D=\Phi(C)$$
is going to play a key role in the following.

We shall estimate degrees, as  the Bogomolov type bound depends on the
degree of the curve.
\begin{lem}
\label{dgr}
$\empty$

\begin{enumerate}
\item
The degree of the curve $\phi(C)$ in $E^2$ is bounded by $6ga^{2}\deg C$.
\item
The degree of the curve $D$ in $E^g$ is bounded by $12g^2a_0^{2(g-2)}a^{2(g-1)}\deg C$.
\end{enumerate}
\end{lem}

\begin{proof}

 i.   Consider   $$\deg \phi(C)=\sum_{i=1}^2\phi(C)\cdot H_i,$$
where $H_i$ is  the coordinate hyperplane given by $3x_i=0$.  The intersection number $\phi(C)\cdot H_i$ is bounded by the degree of the morphism $\phi_{i_{|C}}:C\to E$. Recall that $\phi={\phi_1\choose \phi_2}$.
By Bezout's Theorem, $\deg \phi_{i_{|C}}$ is  at most $3g a^2 \deg C$ - see \cite{V} p. 61.
Therefore $$\deg \phi(C)\le 6g a^2\deg C.$$

ii. Let $X$ be a generic transverse curve in $E^g$.
By    Hindry \cite{Hin} Lemma 6 part i., we deduce
 \begin{equation*}
 \begin{split}
 \deg A^{-1}(X)& \le 2 a^{2(g-2)}\deg X,\\
 \deg A_0^{-1}(X)& \le 2 a_0^{2(g-2)}\deg X.
 \end{split}
 \end{equation*}
 To estimate  $\deg L(X)$, we  proceed  as in  part i.,
 $$\deg L(X)=\sum_{i=1}^gL(X)\cdot H_i,$$ 
where  $H_i$ is  given by $3x_i=0$.  The intersection number $L(X)\cdot H_i$ is bounded by the degree of the morphism $L'_{i_{|X}}:X\to E$, where $L'_i$ is the $i$th. row of $L$.
By Bezout's Theorem $\deg L'_{i_{|X}}$ is  at most $3g a^2 \deg X$ - see \cite{V} p. 61.
Therefore $$\deg L(X)\le 3g^2 a^2\deg X.$$

We conclude that
\begin{equation*}
 \begin{split}
 \deg D\le \deg A_0^{-1}LA^{-1}(C)& \le
  2 a_0^{2(g-2)}\deg LA^{-1}(C)\\
 &\le 6 g^2a_0^{2(g-2)}a^2 \deg A^{-1}(C)\\
 &\le 12 g^2a_0^{2(g-2)}a^2a^{2(g-2)} \deg C.
 \end{split}
 \end{equation*}

\end{proof}
The following Proposition is a lower bound  for  the essential minimum of the image of a curve under Gauss-reduced morphisms.  It reveals the dependence on the height of the morphism.  While the first bound is an immediate application of Conjecture \ref{sin2}, the second estimate is subtle. Our lower bound for $ \mu(\Phi(C+\punto) ) $ grows with  $H(\phi)$. On the contrary, the Bogomolov type lower bound $\epsilon(\Phi(C+\punto) ) $ goes to zero like $\left(a_0H(\phi)\right)^{\frac{-1}{g-1}-\eta}$ - a nice gain.

Potentially, this suggests an interesting question; to investigate the behavior of the essential minimum under a general morphism. 

\begin{propo}
\label{EM}
Assume Conjecture \ref{sin2}. Then,
 for any point $y\in E^g(\overline{\qe})$ and any $\eta>0$, 
 \begin{enumerate}
 \item  
  \begin{equation*}
\mu(\phi(C+\punto))>\ecuno \frac{1}{a^{1+2\eta}},
\end{equation*}
where $\ecuno$ is an effective  constant depending on $C$ and $\eta$. Recall that $a= H(\phi)$.
\item 
  \begin{equation*}
  \mu\left(\Phi(C+\punto)\right )>\ecdue a_0^{\frac{1}{g-1}-\coeta\eta},
  \end{equation*}
  where $\ecdue$ is an effective  constant depending on $C$, $g$ and $\eta$. Recall that $a_0=\lfloor a^\gain\rfloor$.

 \end{enumerate} 

\end{propo}

\begin{proof}
Let us recall  the Bogomolov type bound given in Conjecture \ref{sin2}; for a transverse irreducible curve $X$ in $E^g$ over $\overline{\qe}$ and any $\eta>0$, 
$$\epsilon(X,\eta)=\frac{c(g, E,\eta)}{\deg X^{\frac{1}{2{\rm{cod}}X}+\eta}}< \mu(X).$$

i.  Observe that $\phi(C)\subset E^2$ has codimension $1$.

Let $q'=\phi(\punto)$. So $\phi(C+\punto)= \phi(C)+q'$.  Since $C$ is irreducible, transverse and defined over $\overline{\qe}$,  $\phi(C)+q'$ is as well.
Conjecture \ref{sin2} gives
$$\mu(\phi(C+\punto))= \mu(\phi(C)+q')> \epsilon(\phi(C)+q',\eta)=\frac{c(2, E,\eta)}{\left(\deg (\phi(C)+q')\right)^{\frac{1}{2}+\eta}}.$$

Degrees are preserved by translations,
hence Lemma \ref{dgr} i. implies  that $$\deg (\phi(C)+q')=\deg\phi(C)\le 9ga^2\deg C.$$
If follows that 
$$\epsilon(\phi(C)+q',\eta)\ge \frac{c(2, E,\eta)}{(9ga^2\deg C)^{\frac{1}{2}+\eta}}.$$
Define
$$ \ecuno=\frac{c(2, E,\eta)}{(9g\deg C)^{\frac{1}{2}+\eta}}.$$
Then
$$\mu(\phi(C+\punto))\ge \frac{\ecuno}{a^{1+2\eta}}.$$

ii.
Let $q\in E^g$ be a point such that $[a_0a]q=\Phi(\punto)$. Then
 $$\Phi(C+\punto)=[  a_0a]\left(A_0^{-1}LA^{-1}(C)+q\right)=[a_0a](D+q).$$
 Therefore 
 \begin{equation}
 \label{essmin}
 {\mu\left(\Phi(C+\punto )\right)}=(a_0a)\mu(D+q).
 \end{equation}
 
We now estimate  $\mu(D+q)$ using Conjecture \ref{sin2}.
The curve $D+q$ is irreducible by the definition of $D$. Since $C$ is transverse and defined over $\overline{\qe}$,  $D+q$ is also. 
Thus
$$\mu(D+q)>  \epsilon(D+q,\eta)=\frac{c(g,E,\eta)}{\deg(D+q)^{\frac{1}{2(g-1)}+\eta}}.$$
 Translations by a point preserves degrees, thus Lemma \ref{dgr} ii. gives $$\deg (D+q)=\deg D\le 12g^2a_0^{2(g-2)}a^{2(g-1)}\deg C.$$
 Then 
$$ \epsilon(D+q,\eta)\ge \frac{ c(g,E,\eta)}{(12g^2\deg C)^{\frac{1}{2(g-1)}+\eta}}\left(a_0^{2(g-2)}a^{2(g-1)}\right) ^{-\frac{1}{2(g-1)}-\eta}.$$
Define $$\ecdue= \frac{ c(g,E,\eta)}{(12g^2\deg C)^{\frac{1}{2(g-1)}+\eta}}.$$
So
$$ \mu(D+q)\ge\ecdue{a_0^{-1+\frac{1}{g-1}-2(g-2)\eta}a^{-1-2(g-1)\eta}}.$$
Substitute into (\ref{essmin}), to obtain
$$\mu(\Phi(C+\punto))> \ecdue a_0^{\frac{1}{g-1}-2(g-2)\eta}a^{-2(g-1)\eta}.$$

Recall that $a_0$ is the integral part of $a^\gain$, where $n=2(g+s)-3$. So $2a_0\ge a^\gain$ and
$$a^{2(g-1)\eta}\le (2a_0)^{4n(g-1)\eta}.$$
Further $2(g-2)+4n(g-1) \le 8(g+s)(g-1)$, so
$$\mu(\Phi(C+\punto))> \ecdue a_0^{\frac{1}{g-1}-8(g+s)(g-1)\eta}.$$

\end{proof}

Thankfully we come to our second main proposition; each set in the
union is finite. The proof of i. case (1)  is delicate. In
general $\mu(\pi(C))\le \mu(C)$, for $\pi$ a projection on some factors. We shall rather find a kind of reverse inequality.  On a set of bounded height this will be possible.

\begin{propa*}
Assume Conjecture \ref{sin2}. Then,
there exists $\ecinque>0$ such that 
\begin{enumerate}
\item
For $\varepsilon \le \ecinque$,  for all $\punto
 \in \Gamma_0^2\times \{0\}^{g-2}$  and for
 all Gauss-reduced morphisms $\phi$  of rank $2$, the set 
$$\left(C(\overline{\qe})+\punto\right)\cap \left(B_\phi+\mathcal{O}_{\varepsilon/H(\phi)^\expalt}\right)$$  
is finite.

\item
 For $\varepsilon \le  \frac{\ecinque}{g+s}$ and  for all special morphisms $\tilde{\phi}=(N\phi|\phi')$  of rank $2$, the set 
$$(C(\overline{\qe})\times \gamma) \cap \left(B_{\tilde{\phi}} +\mathcal{O}_{\varepsilon/H(\phi)^\expalt}\right) $$ is finite.

\end{enumerate}

Recall that $n=2(g+s)-3$.

\end{propa*}
\begin{proof}
{\bf Part i.}
Choose 
$$\eta\le \eta_0=\frac{1}{\valeta}.$$
Define  
$$
m=\max \left(2,\left(\frac{\kuno}{\ecdue
  }\right)^{\frac{g-1}{1-8(g+s)(g-1)^2\eta}} \right),$$ 
 and choose
 $$\varepsilon  \le \min \left(\euno,\frac{\kuno}{g}, \frac{ \ecuno }{gm^{\expm}}\right),$$
 where $\ecuno$ and $\ecdue$ are as in Proposition \ref{EM}.  
 
Recall that $H(\phi)=a$. We distinguish two cases:
\begin{itemize}
\item[(1)] $a_0=\lfloor a^\gain\rfloor \ge m$,
\item[(2)]  $a_0=\lfloor a^\gain\rfloor \le m.$
\end{itemize}

Case (1) - \fbox{$a_0\ge m$}

Let $x +\punto  \in  (C(\overline{\qe})+\punto )\cap (B_\phi+\mathcal{O}_{\varepsilon/a^\expalt})$, where $$y=(y_1,y_2,0,\dots,0)\in \Gamma_0^2\times\{0\}^{g-2}.$$ Then
$$\phi(x+y)=\phi(\xi)$$  for  $||\xi||\le \varepsilon/a^\expalt$ and $y\in \Gamma_0^g$.

 We have
chosen $\varepsilon\le \euno$,  so $x \in S_2(C,(\Gamma_0^g)_{\euno})$
which is a set of norm $K_1$. Then \begin{equation*}
||x||\le \kuno.\end{equation*}

Recall that  $\Phi(z_1,\dots,z_g)=(a_0\phi(z),z_3,\dots,z_g) $. So
\begin{equation*}
\begin{split}
\Phi(x+\punto )&=(a_0\phi(x+y),x_3, \dots ,x_g)\\
&=(a_0\phi(\xi), x_3,\dots , x_g).
\end{split}
\end{equation*}
 Therefore
$$ ||\Phi(x+\punto )||=||(a_0\phi(\xi), x_3,\dots , x_g)||\le \max\left(a_0||\phi(\xi)||,||x||\right).$$
Since $||\xi||\le {\varepsilon}{a^{-(\expalt)}}$, $a_0\le a^\gain$ and $\varepsilon \le \frac{K_1}{g}$, we have that  $$a_0||\phi(\xi)||\le a_0(g-r+1)\frac{\varepsilon}{a^\gain}\le K_1.$$
Also $||x||\le K_1$. Thus
 $$||\Phi(x+\punto )||\le \kuno.$$
 
 We work under the hypothesis  $a_0 \ge m= \left(\frac{\kuno}{\ecdue }\right)^{\frac{g-1}{1-8(g+s)(g-1)^2\eta}}$. So 
 \begin{equation*}
  \kuno \le\ecdue a_0^{\frac{1}{g-1}-\coeta\eta}.
   \end{equation*}
 In Proposition \ref{EM} ii., we have proven that
 $$\ecdue
 a_0^{\frac{1}{g-1}-\coeta\eta}<\mu(\Phi(C+\punto)).$$  So  
 $$||\Phi(x+\punto )||\le \kuno< \mu(\Phi(C+\punto)).$$ 
 We deduce that $\Phi(x+\punto )$ belongs to the  finite set $$\Phi(C+\punto)\cap {\mathcal{O}}_{\kuno}.$$ 
 The morphism $C+\punto \to \Phi(C+\punto)$ has finite fiber.
 We can conclude that since $\varepsilon \le\min(\euno,\frac{K_1}{g})$,  for every $\phi$ Gauss-reduced of rank  $2$ with $a_0=\lfloor a^\gain\rfloor\ge m$, the set
$$(C(\overline{\qe})+\punto )\cap \left(B_\phi+\mathcal{O}_{\varepsilon/H(\phi)^\expalt}\right)$$ is finite.

Case (2) -  \fbox{$a_0\le  m$}

 Let $x+\punto  \in (C(\overline{\qe})+\punto )\cap
(B_\phi+\mathcal{O}_{\varepsilon/a^\expalt})$, where $y \in \Gamma_0^2\times\{0\}^{g-2}$.   Then  $$\phi(x+\punto )=\phi(\xi)$$
for $||\xi ||\le \varepsilon/a^\expalt$. However we have chosen  $\varepsilon \le\ecuno /gm^{\expm}$. Hence 
$$||\phi(x+\punto )||=||\phi(\xi)||  \le \frac{g\varepsilon}{a^\gain} \le\frac{ \ecuno }{m^\expm a^\gain}.$$
We are working under the hypothesis 
$a_0=\lfloor a^\gain\rfloor \le m$ and $m\ge 2$, so $a < (2a_0)^{2n}\le m^\expm.$
Furthermore,  $\eta \le \eta_0<\frac{1}{4n}$ implies that $a^{2\eta}<a^\gain$.
Thus
$$a^{1+2\eta}<m^\expm a^\gain.$$
And consequently
$$||\phi(x+\punto )||\le\frac{ \ecuno }{m^\expm a^\gain}<\frac{ \ecuno }{a^{1+2\eta}}.$$
In Proposition \ref{EM} i.  we have proven
$$\frac{ \ecuno }{a^{{1}+2\eta}}<\mu(\phi(C+\punto)).$$
We deduce that  $\phi(x+\punto )$ belongs to the finite set $$\phi(C+\punto)\cap {\mathcal{O}}_{\ecuno m^{-\expm} a^{-\gain}}.$$
The morphism $C+\punto \to \phi(C+\punto)$ has finite fiber.
We conclude that since $\varepsilon \le\frac{\ecuno }{gm^\expm}$,   for all $\phi$ Gauss-reduced of rank  $2$ with $a_0=\lfloor a^\gain\rfloor \le m$, the set 
$$(C(\overline{\qe})+\punto )\cap \left(B_\phi+\mathcal{O}_{\varepsilon/H(\phi)^\expalt}\right)$$ is finite.

For the  curve $C$, define 
\begin{equation*}
\epsilon(C)=\min(\epsilon_1(C,\eta_0),\epsilon_2(C,\eta_0)).
\end{equation*}

Note that 
$\left(\frac{\epsilon(C)
  }{g\kuno}\right)^{8(g+s)g}\le\frac{\ecuno}{gm^{\expm}}$. Thus, we could for instance choose
\begin{equation*}
\label{equattro}
\ecinque=\boundb.
\end{equation*}

\vspace{0.3cm}
 
{\bf Part ii.} 
 We want to show that, for every $\tilde\phi=(N\phi|\phi')$ special of rank $2$, there exists
 $\phi$ Gauss-reduced of
 rank $2$ and $\punto\in \Gamma_0^2\times \{0\}^{g-2}$   such
 that the map $(x,\gamma)\to x+y$ defines an injection
 \begin{equation}
 \label{undici}
 (C(\overline{\qe})\times \gamma) \cap \left(B_{\tilde{\phi}}
 +\mathcal{O}_{\varepsilon/H(\phi)^\expalt}\right) \hookrightarrow 
 (C(\overline{\mathbb{Q}})+\punto )\cap
 \left(B_\phi+\mathcal{O}_{(g+s)\varepsilon/H(\phi)^\expalt}\right). 
 \end{equation}
 We then apply  part i. of this proposition;   if $(g+s)\varepsilon \le \ecinque$, then
$$(C(\overline{\mathbb{Q}})+\punto )\cap
\left(B_\phi+\mathcal{O}_{(g+s)\varepsilon/H(\phi)^\expalt}\right) $$ is finite. So
 if $\varepsilon\le  \frac{\ecinque}{g+s}$, then
$$(C(\overline{\qe})\times \gamma) \cap \left(B_{\tilde{\phi}} +\mathcal{O}_{\varepsilon/H(\phi)^\expalt}\right) $$
is finite.
 
Let us prove the inclusion (\ref{undici}).
Let $\tilde\phi=(N\phi|\phi')$ be special of rank $2$. By definition of special $\phi=(aI_r|L)$ is Gauss-reduced of rank $2$.
Let $$(x,\gamma) \in (C(\overline{\qe})\times \gamma) \cap
\left(B_{\tilde{\phi}} +\mathcal{O}_{\varepsilon/H(\phi)^\expalt}\right) .$$ Then, there exists  $\xi \in \mathcal{O}_{\varepsilon/H(\phi)^\expalt}$ such that  $$\tilde\phi((x,\gamma)+\xi)=0.$$ 
Equivalently $$N\phi(x)+\phi'(\gamma)+\tilde\phi(\xi)=0.$$

 Let ${\punto'}\in E^2$ be a point such that $$N[a]{\punto'}=\phi'(\gamma).$$ 
 Since $\Gamma_0$ is a division group,
 $$\punto=(\punto',0,\cdots ,0)\in \Gamma_0^2\times \{0\}^{g-2}$$ and
 $$N\phi(\punto)=N[a]\punto'=\phi'(\gamma).$$ 
Therefore
$$N\phi(x+y)+\tilde\phi(\xi)=0.$$ 

Let $\xi''\in E^2$ be a point such that $$N[a]\xi''=\tilde\phi(\xi).$$
We define $\xi'=(\xi'', \{0\}^{g-2})$. Then
$$ N\phi(\xi')=N[a]\xi''=\tilde\phi(\xi),$$
and $$ N\phi(x+\punto+\xi')=0.$$
Since $\tilde\phi$ is special $H(\tilde\phi)=Na$. Furthermore $||\xi||\le \frac{\varepsilon}{a^\expalt}$. We deduce
$$||\xi'||=||\xi''||=\frac{||\tilde\phi(\xi)||}{Na}\le
\frac{(g+s)\varepsilon}{a^\expalt}.$$

 In conclusion  $$N\phi(x+\punto +\xi')=0$$ 
 with $||\xi'||\le  \frac{(g+s)\varepsilon}{a^\expalt}$ and $y \in \Gamma_0^2\times \{0\}^{g-2}$.
 Equivalently $$(x+\punto ) \in (C(\overline{\mathbb{Q}})+\punto )\cap \left(B_{N\phi}+\mathcal{O}_{(g+s)\varepsilon/H(\phi)^\expalt}\right) .$$
 By Lemma \ref{tor} i. (with $\psi=\phi$), we deduce 
 $$(x+\punto ) \in (C(\overline{\mathbb{Q}})+\punto )\cap \left(B_{\phi}+\mathcal{O}_{(g+s)\varepsilon/H(\phi)^\expalt}\right) ,$$ 
 with ${\punto}\in \Gamma_0^2\times \{0\}^{g-2}$ and $\phi$ Gauss-reduced of rank $2$.
 
 This proves relation (\ref{undici}) and concludes the proof.

  \end{proof}

\section{The  effectiveness aspect}
\label{effettivo} 

\subsection{An effective weak height bound}

We give an effective bound for the height of 
$S_1(C,\mathcal{O}_\varepsilon)$ for $C$ transverse. 

\begin{thm}
\label{altezza}
Let $C$ be transverse. For every  real $\varepsilon\ge0$, 
the norm of the set
$S_1(C,\mathcal{O}_\varepsilon)$ is bounded by $\kzero\max (1,\varepsilon)$,
 where $\kzero$ is an effective  constant  depending  on  the degree and the height of $C$.
\end{thm}
\begin{proof}
If $x \in S_1(C,\mathcal{O}_\varepsilon)$, there exist $\phi:E^g \to
E$ and $\xi \in  \mathcal{O}_\varepsilon$ such that $\phi (x-\xi)=0$.
Now the proof follows step by step 
 the proof of \cite{V} Theorem 1 page 55 where we replace $\hat{h}$ by $h$, $y$ by $\phi$, $p$ by $x$ and $\hat{h}(y(p))=0$ by $h(\phi(x))=c_0 \deg \phi h(\xi)$ with $h(\xi)\le \varepsilon^2$. 
\end{proof}

\subsection{The strong hypotheses and an  effective weak  theorem }

\begin{proof}[Proof of Theorem \ref{Ctrao}]

The proof  is similar to the proof of Theorem \ref{tutto} given in section \ref{cg}. 

 Theorem  \ref{altezza} implies that for $r\ge 1$ the norm of the set
$S_r(C,\mathcal{O}_1)$ is bounded by an effective constant  $\kzero$.
Let
\begin{enumerate}
\item $\delta_1=\frac{1}{g}\boundaC$.
\item
 $\delta = {\delta_1}{M^{-(\expaltc)}}$ \hfill where $M=\max\left(2,\lceil \frac{\kzero}{\delta_1}\rceil^2\right)^{2g-3}$.\\
\end{enumerate}

In section \ref{subunion}, Proposition $\centro$   i.  with $\Gamma=0$, $\euno=1$ and $\kuno=\kzero$, we have shown that 
 $$\bigcup_{\substack{\phi\,\,\,{\rm{Gauss \,\,\,reduced}}\\ {\rm{rk}}( \phi)=2}}   C(\overline{\qe})\cap \left(B_\phi+\mathcal{O}_{\delta}\right)
 \subset
 \bigcup_{\substack{\phi\,\,\,{\rm{Gauss \,\,\,reduced}}\\ {\rm{rk}}( \phi)=2\,\,\,H(\phi)\le M}}  C(\overline{\qe})\cap \left(B_\phi+\mathcal{O}_{g\delta_1/H(\phi)^\expaltc}\right).$$

 In section \ref{due}, Proposition $\finito$   i. with $\punto=0$ , $s=0$ and $n=2g-3$, $\kuno=\kzero$  and relation (\ref{equattro}), we have shown that for all $\phi$ Gauss-reduced of rank $2$, the set $$ C(\overline{\qe})\cap \left( B_\phi+\mathcal{O}_{g\delta_1/H(\phi)^\expaltc}\right)
$$ is finite. 

The union of finitely many finite sets is finite. It follows that 
$$\bigcup_{\substack{\phi\,\,\,{\rm{Gauss \,\,\,reduced}}\\ {\rm{rk}}( \phi)=2}}  C(\overline{\qe})\cap (B_\phi+\mathcal{O}_{\delta})$$ is finite.

By Lemma \ref{tor1} i. we deduce that
$ S_2(C, \mathcal{O}_{\delta})$ is finite.

This shows that Theorem \ref{Ctrao} holds for
$$\varepsilon\le\boundctra.$$
\end{proof}

\subsection{An effective bound for the cardinality of the sets}

We have  just shown that for $C$ transverse, $\varepsilon$ can be made effective. The purpose of this section is to indicate an effective bound 
   for the cardinality of 
 $S_2(C,\mathcal{O}_{\varepsilon})$, under:
 \begin{con}  [S. David; personal communication]

\label{sin11}

   Let $C$ be a transverse curve in $A$. Then, there exist constants $
   c'=c'(g,\deg_{{{L}}} A, h_L(A), [k:\mathbb{Q}]) $ and $c''=c''(g, \deg_L A, h_L(A), [k:\mathbb{Q}])$
such that, for 
\begin{equation*}
\begin{split}
\epsilon(C) &= \frac{c'}{ (
\deg_{{{L}}} V)^{\frac{1}{2\cod V}}}\\
\Theta(C)&=c''
(\deg_{{{L}}} C)^{g},
\end{split}
\end{equation*}
 the  cardinality of $C(\overline{\mathbb{Q}}) \cap \mathcal{O}_{\epsilon(C)}$ is  bounded by $\Theta(C)$.
\end{con}

This is  the abelian analogue to \cite{fra} Conjecture 1.2.
 
 We prove:
\begin{thm}
\label{ord1} Let $C$ be  transverse.
Assume that Conjecture \ref{sin11} holds.
 Then,  there exists an effective 
$\varepsilon>0$ such that the cardinality of $S_2( C,\mathcal{O}_\epsilon)$
is bounded by an effective constant.
\end{thm}
\begin{proof}
Let $\delta $ and  $\delta_1$  be  as defined  in the previous proof.

By Proposition $\centro$ i.   in section \ref{subunion} we  deduce that
 \begin{equation*}
S_2(C,\mathcal{O}_{\delta})\subset \bigcup_{\substack{\phi \,\,\,{\rm{Gauss-reduced}}\\H(\phi)\le M}}  C(\overline{\qe})\cap \left(B_\phi+\mathcal{O}_{(g+1)\delta_1/H(\phi)^\expaltc}\right).
\end{equation*}
Note that, for any curve $D$ and positive integers $n$, the cardinality of  $
[n]D\cap\mathcal{O}_{n\epsilon(D)}$  is still $\Theta(D)$. Going
through the proofs of Proposition $\finito$  i.  in section \ref{due}, we see that 
$$\sharp S_2(C,\mathcal{O}_{\delta})\le  \sum_{H(\phi)\le M} \sharp\left( \phi^{-1}_{|_C}\left(\phi(C)\cap \mathcal{O}_{\epsilon( \phi(C))}\right)\right),$$
where $\phi_{|_C}:C\to \phi(C)$ is the restriction of $\phi$ to $C$.
Recall that the fiber of $\phi_{|_C}$ has cardinality at most
$3gH(\phi)^2\le 3gM^2$ (see \cite{V} p. 61).
We denote  $${\bf \Delta}_{\rm {max}}=\max_{H(\phi)\le M}\sharp(\phi(C)\cap \mathcal{O}_{\epsilon(\phi(C)}).$$ We deduce 
$$\sharp S_2(C,\mathcal{O}_{\delta})\le 3g M^3{\bf \Delta}_{\rm {max}} .$$
 By Lemma \ref{dgr} i.,$\deg \phi(C) \le (3gH(\phi))^2\deg C$.  Conjecture \ref{sin11} implies that $${\bf \Delta}_{\rm {max}}\le {(3gH(\phi))}^{2g}{\Theta(C)}$$ with $\Theta(C)$ explicitly given.
 We deduce 
 \begin{equation}
 \label{bpcard}\sharp S_2(C,\mathcal{O}_{\delta})\le (3g)^{2g+1}M^{2g+3}{\Theta(C)}.
 \end{equation}
  By Theorem \ref{altezza} the constant $\kzero$ is effective. So $M$ is also  effective. Thus the bound (\ref{bpcard})  is effective, for  $C$  transverse.
 
 \end{proof}
 Similar computations   imply a bound for the cardinality of $S_2(C,\Gamma_{\delta})$.
 
 For $\delta \le \frac{\ecinque}{(g+s)^2}{M'}^{-1-\frac{1}{4g+4s-6}}$,  we  obtain 
 $$\sharp S_2(C,\Gamma_{\delta})\le c_1(g) {M'}^{c_2(g,s)}{\Theta(C)}.
$$
Here $c_1(g)$ (and $c_2(g,s)$)  are  effective constants depending only on $g$ (and $s$). $M'$ depends explicitly on $C$, $g$ and $\kdue$, while $\ecinque$ depends explicitly on $C$, $g$, $s$ and $\kuno$. 

 In view of Theorem \ref{equi34}, the above  bound   also implies a bound for the cardinality of $S_2(C\times \gamma, \mathcal{O}_{\delta/(g+s)\ktre})$.
 
However, Theorem 
\ref{altzero}  does not give    effective    $\kuno$  or $\kdue$. 
Consequently neither $M'$ nor $\ecinque$ are    effective.  An effective estimate  for  $\kuno$ or $\kdue$ would  imply an effective Mordell Conjecture. This gives an indication of  the difficulty to extend effective height proofs  from  transverse curves to weak-transverse curves.

\section{Final Remarks}

\subsection{The C.M.   case}
Sections \ref{nota} - \ref{bound} are proven for any $E$ regardless of whether $E$ has C.M. or not.
Since Conjecture \ref{sin2} is stated for any  $E$, Proposition $\finito$
  holds unchanged for $E$ with C.M.

 We can extend 
 Proposition $\centro$   to Gauss-reduced $\phi \in
M_{r,g}(\mathbb{Z}+\tau \mathbb{Z})$ as follows. Decompose $\phi=\phi_1+\tau\phi_2$ for $\phi_i \in
M_{r,g}(\mathbb{Z})$, then let the morphism $\psi=(\phi_1|\phi_2)$  act on  $(x, \tau x)+(y,\tau y)+(\xi,\tau\xi)$ for $x \in S_r(C,(\Gamma^g_0)_\varepsilon)$, $y \in \Gamma_0^g$ and $\xi \in \mathcal{O}_\varepsilon$.  Apply  Proposition $\centro$ to $\psi$. Constants will depend on $\tau$.

\subsection{From powers to products}
\label{allafine}
 In a power  there are  more algebraic subgroups than in a product where
not all the factors are isogenous. 
If we consider  a product of  non-C.M. elliptic curves, then the matrix of a morphism $\phi$ is simply an integral matrix where the entries corresponding to non isogenous factors are zeros. So nothing changes with respect to our proofs.
If the curve is in a product of elliptic curves in general,  we shall extend the definition of Gauss-reduced, introducing  constants $c_1(\tau)$ and $c_2(\tau)$,  such that the element $a$ on the diagonal has norm satisfying $c_1(\tau)H(\phi)\le |a|\le c_2(\tau)H(\phi)$.

\vspace{0.4cm}

 We leave the details to the reader.

\vskip1cm



\end{document}